\documentclass[11pt,reqno,a4paper]{amsart}

\usepackage[margin=1in]{geometry}
\usepackage{graphicx}
\usepackage{amsaddr}
\usepackage{graphicx, float, changepage, subcaption} 
\usepackage[hidelinks]{hyperref} 
\usepackage[automake=immediate,acronym]{glossaries}
\usepackage{color}
\usepackage{amssymb, amsmath, dsfont,amsthm}
\usepackage{cleveref}
\usepackage{algorithm, algpseudocode}
\usepackage{tikz, pgfplots}
\usepackage{enumitem}
\usepackage{tabls}

\pgfmathdeclarefunction{gauss}{2}{%
	\pgfmathparse{1/(#2*sqrt(2*pi))*exp(-((x-#1)^2)/(2*#2^2))}%
}

\makeglossaries
\newacronym{pde}{PDE}{partial differential equation}
\newacronym{rnla}{RNLA}{randomized numerical linear algebra}
\newacronym{wrt}{w.r.t.}{with respect to}
\newacronym{eg}{e.g.}{for example}
\newacronym{ie}{i.e.}{id est}
\newacronym{mcmc}{MCMC}{Markov Chain Monte Carlo}
\newacronym{eks}{EKS}{Ensemble Kalman Sampler}
\newacronym{cbs}{CBS}{Consensus Based Sampler}
\newacronym{iid}{i.i.d.}{independent and identically distributed}
\newacronym[plural=FIMs,firstplural=Fisher Information Matrices (FIMs)]{fim}{FIM}{Fisher Information Matrix}
\newcommand{\rr}{\mathbb{R}}

\newcommand{\argmin}{\mathrm{argmin}}

\newcommand{\F}{\mathcal{F}}
\newcommand{\C}{\mathcal{C}}
\newcommand{\sfA}{\mathsf{A}}

\newcommand{\sfC}{\mathsf{C}}
\newcommand{\pos}{\mathrm{pos}}
\newcommand{\prior}{\mathrm{pr}}
\newcommand{\eps}{\varepsilon}
\newcommand{\full}{\mathrm{full}}
\newcommand{\inv}{\mathrm{inv}}
\newcommand{\init}{\mathrm{init}}

\newcommand{\EKS}{\mathrm{EKS}}
\newcommand{\CBS}{\mathrm{CBS}}

\renewcommand{\>}{\rangle}

\newcommand{\rd}{\,\text{d}}

\newcommand{\grd}{\nabla}

\newtheorem{theorem}{Theorem}

\newtheorem{remark}{Remark}

\newtheorem{assumption}{Assumption}
\newtheorem{example}{Example}
\makeatletter
\def\namedlabel#1#2{\begingroup
	\def\@currentlabel{#2}%
	\phantomsection\label{#1}\endgroup
}
\makeatother

\title[Sensitivity preserving sampling of the Fisher-Information matrix]{{Local sensitivity-preserving random data down-sampling for experimental design}}

\author{Kathrin Hellmuth$^{a}$, Christian Klingenberg$^{b}$, Qin Li$^{c}$
}
\address{$^{a}$Department of Computing and Mathematical Sciences, California Institute of Technology, USA\\$^{b}$Department of Mathematics, University of W\"urzburg, Germany\\
	$^{c}$Department of Mathematics, University of Wisconsin-Madison, USA}

\email{hellmuth@caltech.edu}
\email{klingen@mathematik.uni-wuerzburg.de}
\email{qinli@math.wisc.edu}

\begin{document}

	\maketitle
	
	\begin{abstract}
	The quality of numerical reconstructions for unknown parameters in inverse problems depends fundamentally on the selection of experimental data. To ensure a robust reconstruction, it is crucial to select data that are sensitive to the parameters, a property typically characterized by the conditioning of the Fisher Information Matrix (FIM). In this work, we propose a general framework for an efficient down-sampling strategy that selects experimental setups that preserve the information content of the full-data FIM. Our approach leverages matrix sketching techniques from randomized numerical linear algebra to achieve a sensitivity-preserving approximation. The method involves drawing samples from a sensitivity-informed distribution, which we execute using gradient-free ensemble sampling methods to handle potentially non-smooth or discrete design spaces. Numerical experiments demonstrate the effectiveness of this framework in selecting optimal sensor locations for a Schr\"odinger potential reconstruction problem.
		
	\end{abstract}
	\begin{quote}
		\noindent 
		{\small {\bf Keywords:} inverse problems, Fisher Information Matrix, Sensitivity analysis, Randomized Numerical Linear Algebra, Sampling methods, Schrödinger potential reconstruction, Ensemble  methods}
	\end{quote}

	\section{Introduction}

	Inverse problems describe the ubiquitous task of inferring a  parameter $p$ from data that follows a statistical model:
	\begin{equation}\label{eqn:measurement}
		y(\xi ) = \mathcal{F}(\xi ,p)\, +\eta(\xi )\,,
	\end{equation}
	where $\mathcal{F}: \Xi \times \mathbb{R}^K\to\mathbb{R}$ {maps the  design variable  $\xi $ and the parameter $p$ to the data. We assume the to-be-inferred parameter $p$ is $K$-dimensional, and $\xi$, the design variable, lives in the underlying design space $ \Xi $ that collects all possible experimental specifications}. The data is polluted by measurement noise {$\eta(\xi)$, which we assume to follow $\mathcal N(0,\Gamma)$ with positive definite covariance matrix $\Gamma\succ 0$}. When the parameter $p$ is fixed, the forward problem returns the solution $y(\xi )$ {for every chosen experimental design variable $\xi $.} 
	
	The associated inverse problem is to revert the process: given evaluations of $y$, we {seek} to infer the parameter $p$.
	There are many approaches to execute this inversion,{ such as cost function minimization, maximum likelihood estimation or Bayesian maximum-a-posteriori estimation. {Denote by $\mathcal R(y,\Xi , \Gamma)$ the chosen reconstruction strategy, and }
		\begin{align}\label{eqn:optimization}
			\hat p = \mathcal R(y,\Xi ,  {\Gamma})
		\end{align}
		{the resulting estimator.}}

	Often there are abundant choices in the design set $\Xi $.  {Assume for simplicity that  $\Xi $ is finite with a large cardinality  $|\Xi | =N\gg K$.}  In this case, it is natural to suspect that one does not need the full data set of $\{y(\xi )\}_{\xi \in\Xi }$ {to characterize $p\in \rr^K$} \cite{alberti2022infinite,harrach2019uniqueness,SimonesPaper}. The task at hand is to select  {a} down-sampled  {data vector $y_c$} that can give an almost equally good recovery of $p$. This reduces experimental as well as computational cost~\cite{chung2019iterative}. More specifically, one aims to design a small finite subset $\Xi _c\subset \Xi $, either through a deterministic or random selection process,   {that} defines the down-sampled data:
	\begin{equation}\label{eqn:down_sample}
		|\Xi _c|=c\ll|\Xi |\,,\quad \text{and}\quad y_c = y|_{\Xi _c}\,,\quad \mathcal{F}_c=\left.\mathcal{F}\right|_{\Xi _c}\,,
	\end{equation}
	{ and endow it with the corresponding  noise model  covariance $\Gamma_c\succ 0$} so that
	\begin{equation}\label{eqn:optimization_c}
		\hat p\quad {\approx}\quad \hat p_{c}=  {\mathcal{R}(y_c,\Xi _c,  {\Gamma_c})}, 
	\end{equation}
	{thereby approximating the recovery problem}~\eqref{eqn:optimization} using a smaller set of data.

	There are many perspectives to take to compare~\eqref{eqn:optimization} and~\eqref{eqn:optimization_c}. {One frequently encountered quantity is  the \gls{fim}. When the full set of data in the design space $\Xi $ is used, the \gls{fim} is defined as:
		\begin{equation}\label{eq:fim:full}
			\begin{aligned}
				\rr^{K\times K} \ni \mathcal{I}(\Xi ) &= \mathcal{I}(\Xi , p_\ast){= G^\top \Gamma^{-1}G}, \\
			\end{aligned}
		\end{equation}
		where $G {= (G_{\xi, :})_{\xi \in \Xi }\in\mathbb R^{N\times K}}$  {with each row collecting} the derivative of $\F(\cdot,\xi)$ at {a given background parameter $p_\ast$}:
		\[
		{G_{\xi, :} }= (\grd_p\F(p_\ast,\xi ))^\top\, \qquad \text{for } \xi \in \Xi .
		\]
		Accordingly, when the data is down-sampled as in~\eqref{eqn:down_sample}, the associated \gls{fim} is given by 
		\begin{equation}\label{eq:fim}
			\begin{aligned}
				\rr^{K\times K} \ni \mathcal{I}(\Xi _c)  {= G_c^\top\Gamma_c^{-1}G_c}, 
			\end{aligned}
		\end{equation}
		where  {$G_c = (G_{\xi, :})_{\xi \in \Xi _c}$}.

		The \gls{fim} is a fundamental tool in optimal experimental design, characterizing the local sensitivity of the data with respect to the parameter vector $p$ in a neighborhood of a nominal parameter $p_\ast$. Through the Cram\'er-Rao inequality (\cite{nordebo2012fisher}), the inverse of the \gls{fim} provides a lower bound on the covariance of any unbiased estimator. In particular, the maximum estimator variance is bounded by the reciprocal of the smallest eigenvalue of $\mathcal{I}(\Xi)$, so larger eigenvalues generally indicate improved local parameter identifiability and reduced estimator uncertainty. Because the \gls{fim} is constructed from a local linearization of the forward map about $p_\ast$, however, it characterizes only the local behavior of the inverse problem. For strongly nonlinear inverse problems, the quality of a design may depend on the global geometry of the parameter-to-observation map, which cannot be fully captured by a single linearization. A common practical strategy is therefore to combine parameter reconstruction and experimental design in a sequential or adaptive manner, repeatedly updating the nominal parameter and recomputing the \gls{fim}, see e.g.~\cite{Santosa}. The present work focuses on understanding and optimizing the linearized inverse problem, which provides the foundation for adaptive strategies in the nonlinear setting.
		
		Once the \gls{fim} has been constructed, two complementary design philosophies can be pursued. The first, and by far the most widely studied, is optimal experimental design (OED), which seeks to optimize a scalar criterion $\Phi$ derived from the spectrum of the \gls{fim}:
		\begin{equation*}
			\Xi_c^{\mathrm{OED}} =\arg \max_{\Xi_c}\Phi(\mathcal{I}(\Xi_c)).
		\end{equation*}
		Prominent examples include E-, A-, and D-optimal designs, which maximize the minimum eigenvalue, the trace, and the determinant of $\mathcal{I}(\Xi)$, respectively \cite{Eswar_Rao_Saibaba,Alexanderian_BayesianOptExpDesign_2021,attia2022optimal,bandara2009optimal,kieferADoptimality1959,mitchell2000algorithm,park2018optimal}. See also the recent review~\cite{huan2024optimalexperimentaldesignformulations}.

		In this work, we take a different perspective, which we refer to as qualitative design \cite{walter1990qualitative,walter1996identifiability}. Rather than seeking an optimal subset of experiments $\Xi_c^{\mathrm{OED}}$ according to a prescribed optimality criterion $\Phi$, our goal is to identify a sufficient subset whose information content preserves the parameter sensitivity of the full experimental configuration $\Xi$. The emphasis is therefore shifted from optimizing $\Phi(\mathcal{I}(\Xi_c))$ to preserving the essential structure of the Fisher information matrix. Mathematically, this amounts to finding those designs $\Xi _c = \{\xi_1,...,\xi_c\}\subset \Xi $ so that $\mathcal{I}(\Xi _c)$ is as informative as $\mathcal{I}(\Xi )$, or
		\begin{center}
			\emph{Design $\Xi _c$ to ensure  {$\mathcal{I}(\Xi _c)\approx \mathcal{I}(\Xi )$}, so that~\eqref{eqn:optimization_c}  holds}.
		\end{center}
		Note that this implies  $\Phi(\mathcal{I}(\Xi_c)) \approx\Phi(\mathcal{I}(\Xi))$ for continuous design criteria.

		In comparison to optimal design where one seeks the best weighting or selection strategy to achieve an optimal eigenvalue structure for the \gls{fim}, this qualitative objective brings us considerable algorithmic flexibility. Specifically, because our goal is to preserve existing sensitivity rather than optimize the spectral profile,  {we can leverage methodologies that are typically underutilized in experimental design}, potentially bypassing the need for computationally expensive iterative solvers. This shift in perspective not only facilitates more efficient computation but also broadens the scope of the theoretical and practical conclusions that can be drawn.

		Indeed, there is a generic condition for $c$, and a generic down-sample strategy that {yields sensitive data} for a very general class of problems. These strategies  are independent of the source of the inverse problem, and do not require any special structure of the original \gls{fim}. The proposed sampling strategy is probabilistic in nature, and thus {sensitivity} is guaranteed with a high probability. This sampling strategy, when applied to any specific problem, leads to a specific distribution for constructing the mask $\Xi_c$. This distribution varies according to the property of $\mathcal{F}$, and thus integrates the knowledge from the underlying model.

		The technical preparation of our approach comes from a seemingly unrelated research area of randomized linear algebra (RNLA)~\cite{mahoney_2016_RandLA,Martinsson_Tropp_2020}. Indeed, the {sensitivity of the data}  is coded in its \gls{fim} \eqref{eq:fim}, {which  enjoys a special tensor structure}. This special structure allows us to deploy random sketching techniques from RNLA to pin down the conditions for preserving the eigenvalue structure. Specifically in this context, we can spell out a probability distribution to draw $\Xi _c$, and show that with high probability, the associated down-sampled {\gls{fim} is well-conditioned,} and thus~\eqref{eqn:optimization_c} holds true with high probability.

		{As such, the novelty of our work  lies in establishing this new perspective on qualitative experimental design in a rather general framework. With this new perspective, we propose a concrete algorithmic pipeline to numerically execute this data selection through sampling.}
		
		Probabilistic methods have been shown to reduce sampling complexity of inverse problems \cite{SimonesPaper} and their integration into  design tasks is an emerging topic that has been studied for instance in \cite{buithan2022bridging,manohar2018datadriven} where the authors used matrix sketching techniques and a low-rank basis representation of the data, respectively. In a Bayesian optimal design setting, a data and model adapted random mask for MRI data acquisition could be constructed in \cite{orozco2024probabilistic} and \cite{Eswar_Rao_Saibaba} applied column subset selection after an initial random down-sampling to reduce computational complexity. Further interesting applications can be found in elliptic solution operator learning \cite{BoulleTownsend_2023_LearnEllipticPDERandLA,boulle2023elliptic} from random input data on the basis of the randomized singular value decomposition. \cite{jin2024uniquereconstructiondiscretizedinverse} examined the same question in a different light, where they used the sketching methods to study how many variables can be stably recovered when the experiments are fixed.

		The two main technical pillars of our proposed method are the matrix sketching, and sampling method. We briefly review them in \Cref{ssec:MatrixSketching} and \Cref{ssec:BayesianSamplingAlgos} respectively. In \Cref{sec:generalProgram} we return to the problems~\eqref{eqn:optimization}-\eqref{eqn:optimization_c}, and examine the relation between their associated \glspl{fim} around the global minimum. The problem will be  {re}cast  { to} invite direct use of random sketching.  {This leads to a very concrete down-sample strategy and provides theoretical guarantees,} as discussed in \Cref{ssec:sampling}. To execute this strategy, practical considerations about sampling choices also play a vital role, and they are discussed in \Cref{ssec:PracticalConsiderations}. In \Cref{sec:applicationtoSchroedinger}, we apply this general program to the potential reconstruction problem for the Schr\"odinger equation, and we conclude the article in \Cref{sec:Discussion}.

		\section{Preview of technical preparations}

		Two main bodies of technical preparation for the current work are matrix sketching techniques rooted in \gls{rnla}, and sampling {algorithms}, rooted in Bayesian problems, that will be leveraged to implement the downsampling. We recall these tools in this section and unify notations.
		
		\subsection{Matrix Sketching in \gls{rnla}}\label{ssec:MatrixSketching}

		\gls{rnla} sees its biggest impact in big data applications, where large data sets, which usually exceed RAM capacities, need to be stored and analyzed quickly. Techniques developed within the domain of \gls{rnla} typically target the access and assessment of a subset of data that is reduced in size but still representative, through ``sketching", see  \cite{huang2021spectral,mahoney_2016_RandLA,Martinsson_Tropp_2020,woodruff2014sketching} and references therein.

		The technique most relevant to our context is the simple computation of matrix-matrix product{s}: how to compute {$\mathsf B:=\mathsf{A}^\top\mathsf{A}\in \rr^{K\times K}$} efficiently? {In the regime where {$\mathsf A {\in \mathbb R^{N\times K}}$ is a} tall but skinny matrix, \gls{ie} when it has significantly  more rows than its number of columns  {$N\gg K$}, it significantly outsizes matrix {$\mathsf B$}, suggesting some row-space information is overly represented.  {A Monte Carlo based method  that sketches  {$\sfA$} by randomly down-sampling rows was proposed in \cite{mahoney_2016_RandLA,Martinsson_Tropp_2020,tropp2025comparison} and is based on the following intuition.}  {The matrix product $\mathsf B$ can be rewritten as a sum}
			\[
			\rr^{K\times K} \ni  {\mathsf B = \mathsf A^\top\mathsf A= \sum_{i = 1}^N \sfA^\top_{i,:} \mathsf A_{i,:}.}
			\]
			{This turns a standard matrix-matrix product into an object written in a summation form, which allows a Monte-Carlo interpretation over a discrete probability measure. This formulation can be relaxed further to allow flexible weights. In particular, choose a user-specified probability distribution $\pi$ on the row indices $\{1,...,N\}$ and define the following random variable}
			\[
			\mathsf{X}= \frac{1}{\pi(i)}\mathsf{A}_{ {i},:}^\top\mathsf{A}_{ {i},:}\,,\quad \text{with }\quad  {i}\sim \pi \,,
			\]
			one can prove: ${\mathsf{B}}=\mathbb{E}\left(\mathsf{X}\right)$.  {Then, following the} standard Monte Carlo  {paradigm, one can replace} the expectation by sample averages:

			\begin{equation}\label{eqn:B_approx}
				{\mathsf{B}}\approx \frac{1}{c}\sum_{j=1}^c\mathsf{X}_j\,,\quad \text{where}\quad \mathsf{X}_j\quad\text{is a {realization of $\mathsf X$}.}
			\end{equation}

			We can summarize this proposal in the following algorithm:
			\begin{algorithm}[H] 
				\caption{BasicMatrixMultiplication from  \cite[Algorithm 3]{mahoney_2016_RandLA}}
				\label{alg:BasicMatrixMult} 
				{\textbf{Input:} $ {N}\times K$ matrix $\sfA$, a  {small} sample size $c\ll N$ and probability measure $\pi$ on $\{1,...,N\}$.\\
					\textbf{Output:} Matrix $\sfC\in \rr^{c\times K}$ such that $\sfC^\top\sfC\approx \sfA^\top\sfA$.}
				\begin{algorithmic}[1]
					\For{$j=1,...,c$}
					\State Sample {$ {i_j}\sim \pi$  i.i.d.; }
					\State Set the $j$-th {row of $\sfC$ as $\sfC_{j:} = \sfA_{ {i_j},:}/\sqrt{c\pi( {i_j})}$.}
					\EndFor \\
					\Return $\sfC$ and $\sfC^\top\sfC$.
				\end{algorithmic}
			\end{algorithm}
			To justify the algorithm,  the approximation sign in~\eqref{eqn:B_approx} can be made more precise, and  the dependence on $c$ and $\pi$ can be spelled out explicitly, largely by deploying central limit theorem and various applications of the Chernoff estimate. It is worth noting that the random variable  $\mathsf{X}$ is a matrix instead of a scalar, so the application of concentration inequality needs caution. Nevertheless, a clever choice of $\pi$ allow{ed the authors in \cite{mahoney_2016_RandLA} to show that  $\sfC^\top\sfC$ approximates $\sfA^\top\sfA$ with high precision and high probability in Frobenius norm:} 
			\begin{theorem}\label{thm:matrixproductclose_mahoney}[\cite[Theorem 7]{mahoney_2016_RandLA}]
				{Let  {$\mathsf A\in \mathbb R^{N\times K}$ with $N\gg K$}. Fix a small, finite number $ {K\leq }c\ll  {N}$ and consider a probability density $\pi$ on  {the row space $\{1,...,N\}$}, for which there exists a $\beta \in (0,1]$ with
					\[
					\pi( {i}) \geq \beta \frac{\|\sfA_{ {i},:}\|_2^2}{\|\sfA\|_F^2}\,,
					\]
					and let the matrix $\sfC\in \rr^{c\times K}$ be constructed by Algorithm \ref{alg:BasicMatrixMult}. Then  for any failure probability rate $\delta \in (0,1)$, one has the error control
				} 
				\begin{align}\label{eq:approxError}
					\mathbb{P}\left(\|\sfA^\top\sfA-\sfC^\top\sfC\|_F\leq \frac{1+\sqrt{8\beta^{-1}\log(\delta^{-1})}}{\sqrt{\beta c}}\|\sfA\|_F^2\right)\geq 1-\delta\,,
				\end{align}
				where the probability $\mathbb{P}$ is taken over all {realizations} of $\sfC$.
			\end{theorem}
			
			The theorem  {studies sampling where} {rows} of $\sfA$  {are chosen from a distribution that is proportional, up to $\beta$,} to its ``volume" -- the $L_2$ norm of the {row}. Then with high probability ($1-\delta$), the approximation of $\mathsf{B}$ by $\sfC^\top\sfC$ is accurate, with the error decaying in the form of $\sqrt{\log{(\delta^{-1})}/c}$  {in Frobenius norm}, where $c$ is the chosen number of sampled rows.
			
			{Note that $\beta\leq 1$, and when $\beta=1$, the mean squared error in~\eqref{eq:approxError} is $\pi(\xi ) = {\|\sfA_{\xi ,:}\|_2^2}/{\|\sfA\|_F^2}$ achieves its minimum.} Suppose we set $\delta = 0.01$. Since $\log(\delta^{-1})=2$  is  $O(1)$, requiring the error to be $\epsilon\|\sfA\|_F^2$ corresponds to $c\geq\frac{O(1)}{\epsilon^2}$. This is an expected MC sampling rate.

			{The concept of expressing a matrix-matrix product as a summation to facilitate Monte Carlo sampling is conceptually simple yet practically powerful. Approximating $\sfA^\top\sfA$ with a sketched version $\sfC^\top\sfC$ aligns directly with our objective of reducing the full \gls{fim} in~\eqref{eq:fim:full} to the compressed form in~\eqref{eq:fim}. However, as shown in~\eqref{eq:fim:full}, the \gls{fim} typically incorporates a noise covariance matrix $\Gamma^{-1}$. It is therefore necessary to extend standard sketching theorems to accommodate the weighted product $\sfA^\top \mathsf{W} \sfA$, particularly when $\mathsf{W}$ is non-diagonal.}

			{The extension is realized by the following observation:
				$$\sfA^\top \mathsf{W} \sfA = \sum_{k,l= 1}^N \sfA_{k,:}^\top \mathsf{W}_{k,l} \sfA_{l,:} = \sum_{k,l= 1}^N \underbrace{\frac{1}{2} \left( \sfA_{k,:}^\top \mathsf{W}_{k,l} \sfA_{l,:} + \sfA_{l,:}^\top \mathsf{W}_{l,k} \sfA_{k,:} \right)}_{\mathsf Y_{k,l}}.$$
				The first equality indicates that, due to the non-diagonality of $\mathsf{W}$, Monte Carlo sampling must be conducted over the double-index $(k,l)$, effectively sketching the row and column spaces. The second equation rewrites the sampling in a way that each component $\mathsf Y_{k,l}$ is symmetric, ensuring that the structural properties of the original \gls{fim} are preserved throughout the sketching process.}
			
			{Following this construction, we define a sequence of i.i.d. random symmetric matrices $\{\mathsf{X}_j\}_{j=1}^c$ according to a user-defined probability distribution $\pi$ over the index space $\{(k,l)\}_{k,l=1}^N$. The weighted product is then approximated as:
				$$\sfA^\top \mathsf{W} \sfA \approx \frac{1}{c} \sum_{j=1}^c \mathsf{X}_j, \quad \text{where} \quad \mathbb{P}\left[\mathsf{X}_j = \frac{1}{\pi(k,l)} \mathsf{Y}_{k,l}\right] = \pi(k,l).$$
				This average can be interpreted as a sketched weighted matrix product, the algorithmic implementation of which is detailed in \Cref{alg:3MatrixMult}.}

			{
				\begin{algorithm}[H] 
					\caption{WeightedMatrixMultiplication}
					\label{alg:3MatrixMult} 
					{\textbf{Input:} $ {N}\times K$ matrix $\sfA$, $ {N}\times N$ matrix $\mathsf W$,  a  {small} sample size $c\ll N$ and probability measure $\pi$ on $\{(k,l)\}_{k,l=1}^N$.\\
						\textbf{Output:} Matrix $\sfC\in \rr^{2c\times K}$ and $\mathsf V\in \mathbb R^{2c\times 2c}$ such that $\sfC^\top\mathsf V\sfC\approx \sfA^\top\mathsf W\sfA$.}
					\begin{algorithmic}[1]
						\State Set $\mathsf V = 0 \in \mathbb R^{2c\times 2c}$
						\For{$j=1,...,c$}
						\State Sample pair {$ {(k_j, l_j)}\sim \pi$  i.i.d.; }
						\State Set the $j$-th and $(j+c)$-th {row of $\sfC$ as $\sfC_{j:} = \sfA_{ {k_j},:}$  and $\sfC_{(j+c):} = \sfA_{ {l_j},:}$.}
						\State Set the respective entries of $\mathsf V$ to $\mathsf V_{j,j+c}=\mathsf V_{j+c,j} = \mathsf W_{k_j, l_j}/(2{c\pi(k_j, l_j)}).$
						\EndFor \\
						\Return $\sfC$, $\mathsf V$ and $\sfC^\top\mathsf V\sfC$.
					\end{algorithmic}
			\end{algorithm}}
			{A direct extension of the arguments in \cite{mahoney_2016_RandLA} suggests that selecting the probability distribution $\pi(k,l)$ proportional to the relative Frobenius norm of $\mathsf Y_{k,l}$ minimizes the mean squared error. This choice enables the derivation of high-probability bounds on the approximation accuracy, as established in the following theorem:
				\begin{theorem}\label{thm:3matrixproductSketch}
					Let $\mathsf W \in \mathbb R^{N\times N}$ be a symmetric positive definite matrix and let $\mathsf A \in \mathbb R^{N\times K}$ with $N \gg K$. For a fixed sample size $c \ll N$, consider a probability distribution $\pi$ over the index space $\{(k,l)\}_{k,l=1}^N$ such that for some $\beta \in (0,1]$, we have:
					$$
					\pi(k,l) \geq \beta \frac{\|\mathsf Y_{k,l}\|_F}{\sum_{k',l'} \|\mathsf Y_{k',l'}\|_F}.
					$$
					Let $\sfC$ and $\mathsf V$ be the sketched matrices constructed according to Algorithm~\ref{alg:3MatrixMult}. Then, for any failure rate $\delta \in (0,1)$, the following high-probability error bound holds:
					\begin{align}\label{eq:approxError_2}
						\mathbb{P}\left( \|\sfA^\top \mathsf W \sfA - \sfC^\top \mathsf V \sfC \|_F \leq \frac{1 + \sqrt{2\beta^{-1} \log(\delta^{-1})}}{\sqrt{\beta c}} \sum_{k,l} \|\mathsf Y_{k,l}\|_F \right) \geq 1 - \delta,
					\end{align}
					where the probability $\mathbb{P}$ is taken over all realizations of the $c$ samples.
				\end{theorem}
				A proof sketch is provided in \Cref{sec:proof3MatrixSketch} for completeness.}

			\subsection{Sampling Algorithms}\label{ssec:BayesianSamplingAlgos}
			
			Sampling  {algorithms are designed to} draw representative samples from a desired distribution (sometimes referred to as target distribution). This task is frequently called upon in the context of Bayesian sampling, where the target distribution is the posterior distribution {with density $\tilde \pi(\xi )\doteq\pi_\pos(\xi |y)\propto l(y|\xi )\pi_\prior(\xi )$ given in terms of the data likelihood $l(y|\xi)$ and the prior density $\pi_\prior$ in finite dimensions.} A {realization} from this distribution provides one {instance of $\xi$ that is likely to have generated the observed data $y$}. In general, due to the positivity of a probability measure, we denote the target distribution
			\begin{equation}\label{eqn:target}
				\tilde \pi(\xi )\propto e^{-\Phi(\xi )}\,,
			\end{equation}
			where $\Phi$ is sometimes referred to as the potential, and $\propto$ means that $\tilde \pi$ is normalized to  integrate to $1$.
			
			Classical methods include \gls{mcmc}-type algorithms. The{ir} strategy is to design a Markov chain whose invariant measure is the target distribution. When a sample walks through this Markov chain, in time, the distribution of the sample converges to the target distribution. Most well-known examples include Langevin Monte Carlo, Hamiltonian Monte Carlo, and Metropolis-Hastings LMC, and so on, compare \cite{bou2013nonasymptotic,chen2014stochastic,cheng2018underdamped,dalalyan2019user,dwivedi2019log,mangoubi2021mixing,RobertsTweedy_1996_ConvergenceMALA}.
			
			Another  paradigm that has recently attracted significant research interest is ensemble-{based}  {sampling}. Originally developed in the context of data assimilation~(\cite{reich2011dynamical,evensen2022data}), this approach has since  {expanded, with}  notable examples  {in sampling}  as summarized in~\cite{Chen2024}, including the \gls{eks}  (\cite{GarbunoInigoHOffmannLiStuart_2020_EKS}) or the \gls{cbs} proposed  in \cite{CarrilloHoffmannStuartVaes_CBS_2022}. These methods evolve an entire ensemble of samples simultaneously through interacting dynamics. The interaction mechanism encodes communication among particles and is carefully crafted to ensure desirable properties, such as being gradient-free or affine-invariant.  Because it remains an  active area of research, non-asymptotic convergence theory  {is} still under development. {It is important to note that all these sampling algorithms do not need the target distribution~\eqref{eqn:target} to be prepared ahead of time. Samples are drawn from an initial guess distribution and updated along the evolution. In the long-time limit, they can be regarded as samples from the target distribution.}
			
			{In the context of Theorem~\ref{thm:3matrixproductSketch}, our task reduces to sampling according to the distribution $\pi(k,l)$. Given the flexibility to choose from a variety of sampling schemes, both classical MCMC and more recent ensemble-based methods are potentially valuable. Because our objective involves identifying a discrete subset of designs $\Xi_c \subset \Xi$, methods that evolve an entire ensemble of particles are particularly relevant. Furthermore, ensemble-based methods typically offer gradient-free formulations. In our setting, the gradient with respect to the design variable $\xi$ can be ill-defined or non-existent (e.g., in discrete design spaces); thus, gradient-free optimization is a highly desirable property. We provide a detailed discussion of the \gls{eks} and \gls{cbs} frameworks below. It should be noted that while these methods were originally developed for continuous state spaces, we employ a nearest-neighbor projection for their numerical realization on our discrete design sets.}
			
			\paragraph{EKS Sampling.}
			\gls{eks} can be viewed as an ensemble version of the Langevin dynamics. It allocates computational resources to update $c$ samples of $\{\xi_j\}_{j=1}^c$ simultaneously:
			\begin{align} \label{EKS:Dynamics}
				\rd{ \xi_j }&=  - C(U)\grd \Phi(\xi_j) \rd{t}+ \sqrt{2C(U)}\rd W_j\,,
			\end{align}
			where $C(U) = c^{-1}\sum_j (\xi_j -\bar{ \xi }) \otimes (\xi_j -\bar{\xi })$ is the empirical covariance matrix between the particles, and $\bar{ \xi } = c^{-1}\sum_{j'} \xi_{j'} $ is the mean.  {The} $W_j$ are independent and identically distributed Brownian motions. Often $\Phi$ takes on a quadratic form: $\Phi(\xi )=\frac{1}{2}\|f(\xi )-y\|^2$. If $f$ is mildly nonlinear,  {then}
			\begin{align}\label{EKS:GradientApprox}
				C(U)\grd \Phi(\xi_j)&=\frac{f(\xi_j)-y}{c}\sum_{j'} (\xi_{j'} -\bar{ \xi }) \otimes (\xi_{j'} -\bar{\xi })\cdot\nabla f(\xi_j)\nonumber\\
				&\approx \frac{f(\xi_j)-y}{c}\sum_{j'} (\xi_{j'} -\bar{ \xi }) \otimes (f(\xi_{j'}) -\bar{f})\,,
			\end{align}
			where the  mild nonlinearity of $f$ allows us to  approximate $\nabla f(\xi_j)$ by a constant for all $\xi_j$. The ensemble average of the $f$ evaluations  {is denoted by $\bar{f}=\frac{1}{c}\sum_jf(\xi_j)$}. Under this weakly nonlinear assumption, the implementation of~\eqref{EKS:Dynamics} is gradient free, and thus achieves a desired property.
			
			When $\Phi$ is Lipschitz-smooth, it was shown in~\cite{ding2021ensemble} and \cite{vaes2024sharp}  that the mean-field limit of~\eqref{EKS:Dynamics} when $c\to\infty$ is:
			\[
			\partial_t\rho=\nabla\cdot(\rho C(\rho)\nabla\Phi) +\mathrm{tr}(C(\rho)D^2\rho\,),
			\]
			and for this equation, it is straightforward to check that $\rho\propto e^{-\Phi}$ is an invariant measure. When $\Phi$ is strongly convex, it was also shown in~\cite{GarbunoInigoHOffmannLiStuart_2020_EKS} that this PDE converges exponentially fast.
			
			In summary,  denoting $ {\rho^c_t}=\frac{1}{c}\sum\delta_{\xi_j {(t)}}$ the empirical distribution, we have, for large enough $c$ and $t$, that  $ {\rho^c_t}\approx \tilde \pi$, and $\{\xi_j\}$ are regarded as samples drawn from the target distribution $\tilde \pi\propto e^{-\Phi}$.

			\paragraph{CBS Sampling.}
			\gls{cbs} was introduced in~\cite{CarrilloHoffmannStuartVaes_CBS_2022} as another method to draw a set of samples simultaneously from a target distribution. It relies on the Laplace principle (\cite{shun1995laplace}). A set of $c$ particles $\{\xi_j\}_{j=1}^c$ evolve according to
			\begin{equation}\label{eq:CBSDynamics}
				\rd \xi_j = - (\xi_j- \mathcal{M}_\beta(\rho_t^c)) \rd t + \sqrt{2 (1+\beta) \Gamma_\beta(\rho_t^c)} \rd W_j\,,
			\end{equation}
			where $\rho^c_t = \frac{1}{c} \sum_{j=1}^c \delta_{\xi_j(t)}$ is the empirical distribution,  {and} $\mathcal{M}_\beta(\rho)$ is the weighted mean parameterized by $\beta$: $\mathcal{M}_{\beta}(\rho) := \mathcal{M}( L_\beta \rho) = \int \xi  \, (  L_\beta \rho)(\!\rd \xi )$.  {Here,}  $ L_\beta \rho= \frac{\rho e^{-\beta \Phi}}{\int\rho e^{-\beta \Phi}\rd \xi }$  {is} the weighted version of $\rho$ and  {the} $\mathcal{M}$ operator takes the mean of a probability distribution. In the $\beta\to\infty$ limit, $ L_\beta \rho$ converges to a Dirac delta centered on the global minimum of $\Phi$ over the support of $\rho$, and thus $\mathcal{M}_{\beta}(\rho)\to\argmin_u\left.\Phi\right|_{\text{supp}(\rho)}$. The second term introduces stochastic deviations in proportion to the covariance of the weighted distribution
			\[
			\Gamma_\beta(\rho) := \Gamma({L}_\beta\rho) := \int \left(\xi -\mathcal{M}(L_\beta\rho)\right)\otimes\left(\xi -\mathcal{M}(L_\beta\rho)\right) \, (L_\beta\rho)(\! \rd \xi )
			\]
			and allows for exploration of the distribution landscape.
			In the mean field limit $c\to\infty$, the particle distribution follows 
			\[
			\partial_t \rho = \grd\cdot \left( \left(\xi -\mathcal{M}(L_\beta\rho)\right) \rho + (1+\beta)  \Gamma_\beta(\rho) \grd \rho\right)\,.
			\]
			Under certain conditions, \cite{CarrilloHoffmannStuartVaes_CBS_2022} showed that the steady state of this equation is a Gaussian approximation of the target distribution around its global maximum, and the PDE solution converges to it exponentially fast. Furthermore, in~\cite{riedl2023gradientneed} the author links this process with Langevin dynamics, viewing it as a gradient-free relaxation.

			\paragraph{ {Early Stopping.}} All of the sampling strategies discussed above can be further improved. One common improvement is to integrate the Metropolis–Hastings (MH) idea as a post-processing. This added step incurs minimal computational cost but helps mitigate bias introduced by the MCMC procedure as described by~\cite{sprungk2023metropolisadjustedinteractingparticlesampling}.

			It is important to note that the introduction of the Metropolis–Hastings (MH) step is primarily aimed at correcting sampling bias. However, other acceptance criteria tailored to the specific problem can also be employed. In our case, for instance, we evaluate the  {conditioning of the down-sampled \gls{fim}}, measured by a selection criterion such as the inverse condition number or the smallest eigenvalue of the  {\gls{fim}}. This simple yet effective strategy is summarized in \Cref{alg:EarlyStoppingSampling} and serves to guide the ensemble evolution toward more favorable configurations through early stopping of the sampling algorithm.

			\begin{algorithm}[H]
				\caption{ {Sampling with early stopping} }
				\label{alg:EarlyStoppingSampling} 
				\textbf{Input:} 
				initial sample $\{\xi_j\}_{j=1,...,c}$, sample update rule $R: {\Xi ^c \to \Xi ^c}$, number of iterations $I>0$, a quantity of interest to be maximized $\mathsf{Q}: {\Xi ^c \to \rr} $\\
				\textbf{Output:} updated sample $\{\xi_j\}_{j=1,...,c}$ with improved evaluation criterion.
				\begin{algorithmic}[1]
					\For{$i=1,...,I$}
					\State Generate sample update: $\{\theta_j\}_{j=1,...,c} = R(\{\xi_j\}_{j=1,...,c} ).$
					\If{$\mathsf{Q}(\{\theta_j\})>\mathsf{Q}(\{\xi_j\})$, } 
					Update $\{\xi_j\}_{j=1,...,c}\leftarrow \{\theta_j\}_{j=1,...,c}$
					\EndIf
					\EndFor \\
					\Return sample $\{\xi_j\}_{j=1,...,c}$.
				\end{algorithmic}
			\end{algorithm}

			\section{The general program} \label{sec:generalProgram}
			Having reviewed both matrix sketching techniques and sampling algorithms, we now return to our qualitative experimental design problem and apply these tools to address it. Specifically, our goal is to identify suitable experimental setups that preserve data sensitivity for parameter reconstruction~\eqref{eqn:optimization}, even when the data is down-sampled~\eqref{eqn:optimization_c}. To achieve this, we reformulate the task as a sketching problem over the {\gls{fim}}, allowing us to leverage results such as~\Cref{thm:3matrixproductSketch} for theoretical guarantees. This reformulated problem, when executed numerically, is coupled with a sampling strategy. In particular, ensemble-based sampling methods - such as those described in~\eqref{EKS:Dynamics} and~\eqref{eq:CBSDynamics} - are employed to guide the selection process.
			
			In the following sections, we begin by revisiting the structure of the \gls{fim}, which lays the foundation for applying sketching techniques~\Cref{ssec:setups}. This is followed by the integration of sampling methods as an algorithmic strategy for selecting informative data points~\Cref{ssec:sampling}. Additional practical considerations are discussed in~\Cref{ssec:PracticalConsiderations}.

			\subsection{Setup and Sampling Perspective}\label{ssec:setups}
			{Recall that the full-data \gls{fim} and the subsampled \gls{fim} are defined as:
				\[
				\mathcal{I}(\Xi, p_\ast) = G^\top \Gamma^{-1}G \quad \text{and} \quad \mathcal{I}(\Xi_c) = G_c^\top \Gamma_c^{-1} G_c\,,
				\]
				where our objective is to identify a subset $\Xi_c$ such that $\mathcal{I}(\Xi_c)$ accurately approximates $\mathcal{I}(\Xi, p_\ast)$. To quantify this approximation, we denote the minimum eigenvalues of $\mathcal{I}(\Xi)$ and $\mathcal{I}(\Xi_c)$ as $\lambda_{\min}^\Xi$ and $\lambda_{\min}^c$, respectively, and their corresponding inverse condition numbers as $c_{\text{inv}}^\Xi$ and $c_{\text{inv}}^c$. A down-sampled $\Xi_c$ is considered successful if $\lambda_{\min}^c$ is strictly positive and $c_{\text{inv}}^c$ is close to $1$.}
		}
		
		{To proceed, we establish a primary assumption that defines our operational framework:}
		
		\begin{assumption}\label{ass:groundtruthwellposed}
			{The experimental design $\Xi$ is locally sensitive at $p_\ast$ in the sense that the \gls{fim} defined in~\eqref{eq:fim:full} is of full rank for a given a positive-definite covariance $\Gamma$. The inverse problem is therefore theoretically well-posed.}
		\end{assumption}
		
		This assumption is foundational; it ensures that the full dataset contains sufficient information to allow for a unique parameter reconstruction, providing a benchmark for our subsampling strategy.
		
		This formulation aligns with the matrix sketching framework described in~\Cref{ssec:MatrixSketching}. Specifically, preserving parameter sensitivity under design down-sampling translates to selecting rows from the sensitivity matrix $G := G(p_\ast)$ and constructing a noise precision matrix $P$ such that the symmetric approximation
		\begin{equation}\label{eqn:G_star_approx}
			{\sum_{\xi, \zeta \in \Xi} G_{\xi, :}^\top \Gamma_{\xi, \zeta}^{-1} G_{\zeta, :} \approx \frac{1}{2c} \sum_{\xi, \zeta \in \Xi_c} (G_{\xi, :}^\top P_{\xi, \zeta} G_{\zeta, :} + G_{\zeta, :}^\top P_{\zeta, \xi} G_{\xi, :})}
		\end{equation}
		holds with high probability.  {In this context and throughout the following sections, we identify the rows of $G$ with their corresponding experimental specifications, i.e., $G_{i, :} = G_{\xi_i, :}$, for a fixed ordering of the design space $\Xi = \{\xi_1, \dots, \xi_N\}$.}

		\subsection{Experimental Design through Sampling} \label{ssec:sampling}
		In light of~\eqref{eqn:G_star_approx}, we deploy~ {\Cref{alg:3MatrixMult}} and obtain the following down-sampling strategy:
		{
			\begin{align}\label{Ccw:def}
				\text{for }(\xi_j,\theta_j)\sim \pi:\quad  
				&\Xi _c=\{\xi_j\}_{j=1}^c \cup\{\theta_j\}_{j=1}^c\quad\text{with}\quad |\Xi _c| {\leq 2c}, \qquad 
				y_c= y|_{\Xi _c}
				\,,\\
				&   P\in \mathbb R^{2c\times 2c} \quad \text{with}\quad  P_{i,j}= \begin{cases}\frac {\Gamma^{-1}_{(\xi_j,\theta_j)}} {c\pi(\xi_j,\theta_j)}, &\text{for } |j- i|=c \\
					0, & \text{else.}
				\end{cases}  .\nonumber
			\end{align}
		}

		The associated {\gls{fim}} at the  {nominal parameter} then becomes:
		\begin{equation}\label{eq:fimweight}
			\tilde {\mathcal I} (\Xi _c) =   {G_c^\top PG_c = \sum_{\xi ,\theta \in\Xi _c} \frac{1}{2c\pi(\xi ,\theta )}\underbrace{\left(G_{\xi ,:}^\top \Gamma^{-1}_{\xi ,\theta }G_{\theta ,:}+G_{\theta ,:}^\top \Gamma^{-1}_{\theta ,\xi }G_{\xi ,:}\right)}_{=: \mathsf Y_{\xi,\theta}\,.}}
		\end{equation}

		As suggested in  {\Cref{thm:3matrixproductSketch}}, there is an optimal sampling strategy with each  { pair of experiments $(\xi ,\theta )$} being selected with a rate proportional to  {the} volume  {of its contribution to $\mathcal I(\Xi )$}. In our context, this optimal strategy is:
		\begin{equation}\label{eq:ptilde}
			{\tilde \pi(\xi ,\theta )\propto\|
				\mathsf Y_{\xi, \theta}
				\|_F\,.}
		\end{equation}

		If any chosen $\pi$ is close to the optimal $\tilde\pi$, {\Cref{thm:3matrixproductSketch} ensures that with high probability, the down-sampled data is locally sensitive to the parameter}  if the sample size $c$ is sufficiently large. {The following theorem provides a probabilistic bound to preserve data sensitivity.}
		
		\begin{theorem}\label{thm:SamplingDesigns}
			Consider an inverse problem that satisfies Assumption \ref{ass:groundtruthwellposed},  and let the re-weighted {data}  be constructed  as in \eqref{Ccw:def}, where the sampling probability density $\pi(\xi , {\theta })$ on $\Xi  {\times \Xi }$ satisfies 
			\begin{align}\label{pi_n:bound}
				\pi \geq \beta \tilde \pi,
			\end{align}
			for some $\beta \in (0,1]$. Then, for any failure probability $\delta \in (0,1)$ and any error tolerance $\eps \in \left(0,\lambda_{\min}^\Xi \right)$, a choice of the sample size
			\begin{equation}\label{eq:cBound}
				c\geq  {\bigg(\sum_{\xi,\theta\in \Xi}\|Y_{\xi,\theta}\|_F\bigg)^2}\frac{(1+\sqrt{ {2}\beta^{-1}\log(\delta^{-1})})^2}{\beta {\eps^2}}
			\end{equation}
			assures that:
			\[
			\lambda_{\min}^c\geq {\lambda_{\min}^\Xi -} \eps>0\,,\quad\text{and}\quad c_{\inv}^c\geq c_\inv^\Xi  \frac{\lambda_{\min}^\Xi - \eps}{{\lambda_{\min}^\Xi }+\eps}\,.
			\]
		\end{theorem}
		
		\begin{proof}
			Noting that Theorem~\ref{thm:3matrixproductSketch} gives the control over Frobenius norm, and thus the spectral norm is a direct consequence. In particular:
			\begin{align}\label{eq:lambdaminHCc}
				\lambda_{\min}^c \geq& \lambda_{\min}^\Xi  - |\lambda_{\min}^\Xi -\lambda_{\min}^c|  
				\geq  \lambda_{\min}^\Xi  - \|\mathcal I (\Xi ) - \mathcal I (\Xi _c)\|_F\,,
			\end{align}
			where the second term is controlled by~ {\Cref{thm:3matrixproductSketch}}. Namely, with probability at least $1-\delta$
			\begin{align*}
				\|\mathcal I (\Xi )- \mathcal I (\Xi _c)\|_F &
				\leq \frac{1+\sqrt{ {2}\beta^{-1}\log(\delta^{-1})}}{\sqrt{\beta c}} {\sum_{\xi,\theta\in \Xi}\|\mathsf Y_{\xi,\theta}\|_F}\,.
			\end{align*}
			Setting it to be {$\eps$}, we obtain~\eqref{eq:cBound}. Similar estimates of the maximum eigenvalue yields the bound for the inverse condition number $c_\inv^c = \frac{\lambda_{\min}^c}{\lambda_{\max}^c}.$
		\end{proof}
		
		{\begin{remark}[Accessibility of $\tilde{\pi}$] \label{rem:accessPiTilde}
				Evaluating the optimal sampling strategy $\tilde{\pi}(\xi, \theta)$ in \eqref{eq:ptilde} requires access to the following components:
				\begin{enumerate}[label = (\roman*)]
					\item The Inverse Covariance Matrix $\Gamma^{-1}$: Computing $\Gamma^{-1}$ can be computationally prohibitive if the noise model is defined directly via the covariance matrix $\Gamma$. However, in many practical settings, the noise model is explicitly prescribed through the precision matrix (the inverse covariance), making its entries directly accessible.
					\item The Sensitivity Matrix $G$: It measures sensitivity of the data locally around the given background parameter $p_\ast$, which  is  meaningful for reconstruction only when $p_\ast \approx p_0$ coincides, at least approximately, with the data-generating ground-truth parameter $p_0$.     This requirement highlights a classical challenge in optimal experimental design: the optimal selection of measurements depends on the very parameter one seeks to estimate. In practice, this often necessitates an iterative or ``online" approach, where the parameter estimate and the experimental design are updated alternately. While such greedy-type updates are essential for non-linear problems, they remain beyond the scope of the current work; we refer the reader to~\cite{huan2024optimalexperimentaldesignformulations,Alexanderian_BayesianOptExpDesign_2021} for comprehensive reviews on these adaptive strategies.
				\end{enumerate}
			\end{remark}
		}
		
		\begin{example}\label{example:uncorrelated}
			{When noise is uncorrelated,  the covariance $\Gamma_{ij} =\gamma_{ii}\delta_{ij}$ is diagonal, and thus $\mathsf Y_{\xi ,\theta }\neq 0$ only if $\xi = \theta$. This means the support of $\tilde\pi$ is only on the trace of $\xi=\theta$. The method is simplified, with the sampling probability becoming
				\begin{equation}\label{eq:uncorrelatedSketching}
					\tilde \pi(\xi ) \propto \|G^\top_{\xi ,:}\Gamma^{-1}_{\xi ,\xi }G_{\xi ,:}\|_F = \|\Gamma^{-1/2}_{\xi ,\xi }G_{\xi ,:}\|_2^2\,.
				\end{equation}
				This is in line with row sketching of the reweighted gradient matrix $J:= \Gamma^{-1/2} G$ for design with the classical algorithm \Cref{alg:BasicMatrixMult}.}
		\end{example}
		\begin{remark}[Sampling infinite  designs]
			{In practice, many design spaces $\Xi $ are continuum-indexed. For instance, sensor locations can be anywhere in a compact domain. Our results can be generalized to these infinite-dimensional settings with $G$ replaced by a quasi-matrix. In situations like that, one needs to prescribe the underlying measure $\mu(\xi)$ for the design space.}
		\end{remark}
		
		{
			\begin{algorithm}[H] 
				\caption{Experimental Design sampling}
				\label{alg:SampleDesigns} 
				{\textbf{Input:} probability distribution $\pi$ on design ind{ex pairs} $\{(k,l)\}_{k,l=1}^N$, desired number ${2c}\ll N$ of experiments, {noise precision matrix $\Gamma^{-1}$}.\\
					\textbf{Output:} Down-sampled designs {$\{\xi_1,...,\xi_c, \theta_1,...,\theta_c\}$.}}
				\begin{algorithmic}[1]
					\State {Set $P \leftarrow 0\in \mathbb R^{2c\times 2c}$}
					\For{$j=1,...,c$}
					\State Sample {${(\xi_j,\theta_j)}\sim \pi$  i.i.d.; }
					\State Add {indices $\xi_j,\theta_j$} to set of indices $\Xi _c$;
					\State {Store new precision value $P_{j,j+c}= P_{j+c,j} \leftarrow {\Gamma^{-1}_{\xi_j,\theta_j}}/({c\pi(\xi_j,\theta_j)})$.}
					\EndFor \\
					\Return Down-sampled experimental designs $\{\xi_j,\theta_j\}_{j =1,...,c}$ and  noise precision matrix $P$.
				\end{algorithmic}
			\end{algorithm}
		}
		
		\subsection{Practical considerations}\label{ssec:PracticalConsiderations}

		According to~\Cref{thm:SamplingDesigns}, we are  {seeking to draw}  $c$ i.i.d. samples from the  {sensitivity-based sampling}  distribution $\tilde \pi$.  {While in principle, both MCMC-type and ensemble-type sampling algorithms can be applied, ensemble-type sampling methods usually have their gradient-free versions and are thus more suitable for our setting. Here, we adopt the algorithms described in~\Cref{ssec:BayesianSamplingAlgos}.}
		
		{
			Assume that $\tilde \Xi $ is of finite Lebesgue measure and that $G_{\xi, :}$ and $\Gamma^{-1}_{\xi ,\theta }$ are  uniformly bounded for all $\xi ,\theta \in \tilde \Xi ^2$ and   non-zero. Then  $\|\mathsf Y_{\xi ,\theta }\|_F$ is bounded and non-zero. Our target distribution is}
		{
			\begin{align}\label{phi:def}
				\tilde \pi(\xi ,\theta )= \frac{1}{Z}e^{-\Phi(\xi ,\theta )} \quad \text{ with } \quad  \Phi(\xi ,\theta ):= -\log\left(\|\mathsf Y_{\xi ,\theta }\|_F\right).
			\end{align}
			where $Z$ is the normalization constant.
		}
		
		A direct application of EKS provides us the following  strategy:  {Start with a simple distribution e.g., uniform or Gaussian, over $\tilde\Xi$ and draw $c$ samples from it to form ${S}= \{\chi_j = (\xi_{j},\theta_j)\}_{j=1,...,c}$}. These $c$ samples are then evolved interactively according to:
		\begin{align}\label{eq:EKS_samples} 
			{ \rd \chi_j =  \sum_{k}  D_{j,k} \chi_{k} \rd t+ \sqrt{2C(S)}\rd W_j}\,,
		\end{align}
		{where the interaction kernel incorporates information from $\tilde \pi$, and is defined as}
		{
			\begin{align*}
				D_{j,k} :=&
				\frac{1}{c\|\mathsf Y_{\chi_j}\|^2_F}\langle \mathsf Y_{\chi_{k}} - \overline{ \mathsf Y_{\chi}}, \mathsf Y_{\chi_j}\rangle_F
			\end{align*}
			and $C(S) = \frac{1}{c}\sum_{k}(\chi_{k} - \bar \chi)\otimes(\chi_{k} - \bar \chi)$ is the covariance matrix.
		}
		{Note that this term is approximately
			\begin{align*}
				\sum_{k}  D_{j,k} \chi_{k}
				\approx&\left(\frac{1}{2c\|\mathsf Y_{\chi_j}\|^2_F}\sum_{k}(\chi_{k} - \bar \chi)\otimes(\chi_{k} - \bar \chi)\right)\cdot \nabla_\chi \langle \mathsf Y_{\chi_{j}}, \mathsf Y_{\chi_j}\rangle_F \\
				=&
				-\left(\frac{1}{c}\sum_{k}(\chi_{k} - \bar \chi)\otimes(\chi_{k} - \bar \chi)\right)\cdot\nabla_\chi\Phi(\chi_j)=-C(S)\cdot\grd_{\chi} \Phi(\chi_j).
			\end{align*}
			and approximates the gradient information, where we used a derivation in analogy to that in \eqref{EKS:GradientApprox}. So the invariant measure of the SDE recovers the target distribution~\eqref{eq:EKS_samples}.
		}
		Running this SDE forward in time using the classical Euler-Maruyama method, we get the updated samples:
		{
			\begin{align}\label{eq:EKS_Timestep}
				\chi_j^{t_{n+1}} =&\ \chi_j^{t_n} + \Delta t_n \sum_{j'} D_{j,j'}^{t_n} \chi_{j'}^{t_n}
				+ \sqrt{2\Delta t_n C(S^{t_n})} \eta_j^{t_n}\,,
			\end{align}
		}
		with $ {\eta_j^{t_n}} \sim N(0,I)$ independent and identically distributed. For  enhanced stability, numerically one can further deploy adaptive schemes proposed in \cite{Kovachki_2019,GarbunoInigoHOffmannLiStuart_2020_EKS} by setting $\Delta t_n= \frac{\Delta t_0}{\|D^{t_n}\|_F+ \eps}$ for some $\eps>0$ where $D^{t_n}= (D^{t_n}_{j,k})_{j,k}$.

		Application of \gls{cbs} is equally straightforward. Following~\cite{CarrilloHoffmannStuartVaes_CBS_2022}, we deploy the forward in time discretization using an exponential integrator:
		\begin{equation}\label{eq:CBS_Timestep}
			\chi_j^{t_{n+1}} = e^{-\Delta t} \chi_j^{t_n} + (1-e^{-\Delta t} )\mathcal M_\beta(\rho^{c}_{t_n}) + \sqrt{(1-e^{-2\Delta t}) (1+\beta)\Gamma_\beta(\rho^{c}_{t_n})} \eta_j^{t_n}.
		\end{equation}
		{noting that sensitivity information from $\tilde \pi$ enters through the reweighting of the empirical distribution $\rho^c_t = \frac{1}{c} \sum_{j=1}^c \delta_{\xi_j(t)}$ via  $ L_\beta \rho= \frac{\rho e^{-\beta \Phi}}{\int\rho e^{-\beta \Phi}\rd \xi }$ in the computation of the reweighted mean and covariance as described underneath \eqref{eq:CBSDynamics}:
			\begin{align*}
				&\mathcal{M}_{\beta}(\rho) := 
				\int \xi  \, (  L_\beta \rho)(\!\rd \xi ) \qquad \text{and}\\
				& \Gamma_\beta(\rho) 
				:= \int \left(\xi -\mathcal{M}_\beta(\rho)\right)\otimes\left(\xi -\mathcal{M}_\beta(\rho)\right) \, (L_\beta\rho)(\! \rd \xi ).
			\end{align*}
	}}

	\begin{remark}
		{In the case when $\Gamma$ is diagonal, noise is uncorrelated. According to Example~\ref{example:uncorrelated}, the optimal distribution takes  the form:
			\begin{align}\label{phi:def}
				\tilde \pi(\xi )= \frac{1}{Z}e^{-\Phi(\xi )} \quad \text{ with } \quad  \Phi(\xi ):= -\log\left(\|J_{\xi ,:}\|_2^2\right)\,\qquad \text{for } J= \Gamma^{-1/2} G.
			\end{align}
			In this case, the sampling of $\chi$ is reduced to that of $\xi$, and  the applications of both, EKS and CBS, become simplified.}
		
		{Take EKS for example, upon preparing $c$ interactive samples $U= \{\xi_{j}\}_{j=1,...,c}$ from an initial sample distribution such as uniform or Gaussian, the reduced EKS~\eqref{eq:EKS_samples} becomes $\rd \xi_j =  \sum_{k}  D_{j,k} \xi_{k} \rd t+ \sqrt{2C(U)}\rd W_j$, where $C(U)$ is the empirical covariance matrix and $D$ encodes $\tilde\pi$ information and is defined as $D_{j,k}= \frac{1}{c\|J_{\xi_j,:}\|_2^2}\<J_{\xi_{k,:}}- \overline{J_{\xi,:}}\,, J_{\xi_j,:}\>_2 $.}
	\end{remark}

	\subsection{Drawbacks of the Sampling Method and Early Stopping}\label{ssec:earlystopping}
	{While the proposed sampling methods are effective, several practical drawbacks warrant discussion.}
	
	{First, many ensemble-based methods are still in the early stages of theoretical development. The absence of non-asymptotic convergence rates, particularly when the target distribution is non-Gaussian, presents a challenge for performance guarantees. However, the gradient-free nature of these methods offers a significant computational advantage in scenarios where the gradient of the forward model is either unavailable or prohibitively expensive to compute.}
	
	{Second, the target distribution $\tilde{\pi}$ depends on quantities such as $\sum_{\xi, \theta \in \Xi} \|\mathsf Y_{\xi, \theta}\|_F$ and $\lambda_{\min}^\Xi$, which are typically unknown a priori and must be estimated dynamically. While this complicates the sampling process, Theorem \ref{thm:SamplingDesigns} establishes that ``perfect" sampling is not required; the framework is robust to a discrepancy of $\beta < 1$. Under the heuristic assumption that $G_{\xi, :}$ and $\Gamma^{-1}_{\xi, \theta}$ remain uniformly bounded, a uniform density serves as a sufficient proxy for $\tilde{\pi}$. This robustness is corroborated by our numerical results in \Cref{fig:Sampling_UniformData}.}

	{Finally, it is worth emphasizing that our objective differs fundamentally from classical Bayesian sampling. Rather than seeking high-fidelity posterior samples, our goal is to improve the conditioning and ensure the positivity of the \gls{fim}. Because sampling accuracy is secondary to numerical stability, we utilize an early stopping criterion: the sampling procedure is terminated once a favorable condition number is achieved. This heuristic significantly reduces the overall computational overhead, as evidenced by the evolution of $c_{\mathrm{inv}}$ in \Cref{fig:evolveLambdaminEvolution_NormalData,fig:SamplingStatistics}.}
	
	{We summarize the procedure in~\Cref{alg:SampleDesigns_EnsembleGreedy}.}

	{
		\begin{algorithm}[H] 
			\caption{Experimental Design sampling with ensemble methods and  early stopping}
			\label{alg:SampleDesigns_EnsembleGreedy} 
			{\textbf{Input:} Initial data distribution $\pi_0$, algorithm to compute $G_{\xi, :}$ for given $\xi $, desired number $c\ll N$ of experiments, number of iterations $I>0$, evaluation criterion: scalarization of $\mathcal I (\Xi _c)$ to be maximized $q:\rr^{K\times K} \to [0,\infty)$.\\
				\textbf{Output:} updated sample  $\{\xi_j\}_{j=1,...,c}$ with an improved evaluation criterion.}
			\begin{algorithmic}[1]
				\State Initialize ensemble $\Xi _c = \{\theta_j\}_{j=1,...,c}$ with $\theta_j\sim \pi_0$, $j=1,...,c$;
				\State Initialize best seen ensemble $\Xi _c^\ast \leftarrow \Xi _c$
				\State Compute $G_c = (G_{\theta_j,:})_{\theta_j\in \Xi _c}$; 
				\State Build $\tilde{\mathcal I}(\Xi _c)$ according to   \eqref{eq:fimweight} with $\pi\mid_{\Xi _c} = \pi_0\mid_{\Xi _c}$;
				\State Initialize best evaluation criterion $q^\ast \leftarrow q(\tilde{\mathcal I}(\Xi _c))$
				\For{$i=1,...,I$} 
				\State Update $\Xi _c$, either by \eqref{eq:EKS_Timestep} or \eqref{eq:CBS_Timestep}, using $G_c$;
				\State Update $G_c = (G_{\theta_j,:})_{\theta_j\in \Xi _c}$ on the new ensemble; 
				\State Compute $\tilde {\mathcal I}(\Xi _c)$ according to  \eqref{eq:fimweight}
				\If{$q(\tilde {\mathcal I} (\Xi _c) )>q^\ast$} Update  $\Xi _c^\ast \leftarrow \Xi _c$ and  $q^\ast \leftarrow  q(\tilde {\mathcal I} (\Xi _c) )$;
				\EndIf
				\EndFor \\
				\Return Best seen ensemble $\Xi _c^\ast$.
			\end{algorithmic}
		\end{algorithm}
	}

	\section{Application to   Schr\"odinger potential reconstruction} 
	\label{sec:applicationtoSchroedinger}
	In this section, we numerically examine the performance of the proposed algorithm using the steady-state Schr\"odinger equation as a specific example. The spatial domain is set to be $X=[-1,1]^2\subset\rr^2$ and the PDE is equipped with non-negative source term $\gamma\in C^\infty_+(X)$:
	\begin{align}\label{forward}
		(-\Delta +  p)u_p = \gamma &\quad x\in X,\\
		u_p = 0&\quad x\in\partial X.\nonumber
	\end{align}
	{ {Let $\gamma\neq 0$. Then} the existence of a positive solution $u_p$ for a fixed non-negative parameter $p\in C^\infty_{+}(X)$ follows from standard elliptic theory. }
	
	The inverse problem is to reconstruct the potential $p$ from measurements of the observable solution $u_p$. Clearly inferring $p$ becomes trivial when one has the full noise free data $u_p$: $p=\frac{\gamma+\Delta u_p}{u_p}$ pointwise in $X$. Problem{s arise}  in the finite dimensional setting  {where o}nly a finite number of {potentially noisy} measurements of $u_p$ is taken and $p$ is parameterized by finite{ly} many parameters. The goal is to find  {a suitable}  experimental setting (measuring location) for  inferring $p$.

	\paragraph{Parameter Discretization.}
	Let $\{\phi_k:X\to \rr\}_{k=1,...,K}$ be a  given finite set of  {smooth} basis functions on $X$. Our admissible set {consists of non-negative functions $p$ of the form}:
	\begin{align*}
		\mathcal{A}:= \bigg\{p:X\to \rr_0^{+}, \quad  x = \begin{pmatrix}
			x_1\\x_2
		\end{pmatrix} \mapsto p(x) &= \sum_{k=1}^K p_k \phi_k(x_1, x_2)\\
		&\text{ for some }p_k\in \rr, k=1,...,K\bigg\}.
	\end{align*}
	In the numerical examples below, we used $K=9$ with basis 
	\[
	\{\phi_{k_1,k_2}(x_1, x_2) = \cos(k_1 \pi x_1) \cos(k_2\pi  x_2)\}_{k_1,k_2= 0,1,2}\,.
	\]

	\paragraph{Space discretization.} To numerically realize the PDE solution, we use its numerical solution computed on  {an} equidistant Cartesian grid $\{x_n, n = 1,...,(N_x+1)^2\}$, where we set $N_x$ cells in every direction.

	\subsection{Fixed Source Term $\gamma$}
	In the first set of experiments, we fix the source term at $\gamma = 10^4$ for all trials. Based on the above considerations, this well-controlled setting should be sufficient for successful reconstruction, provided that the data is  chosen appropriately.

	\paragraph{Data.} 
	We assume all possible measurements are point-wise measurements, meaning $\F(x, p)=u_p(x)$ for all $x\in X$. We denote by  $p_0 \in \mathcal{A}$ the  ground-truth potential   and  assume {that the noise} satisfies $\eta(x)\sim\mathcal N(0,1)$ \gls{iid}. The resulting data are
	\[
	\{y(x)= \F(x,p_0)+ \eta(x)=u_{p_0}(x)+\eta(x)\}_{x\in \Xi}\,.
	\]
	The question related to experimental design now translates into  finding  the number $c$ and locations $\Xi_c\subset X$ such that the associated down-sampled optimization problem is locally strictly convex. As a consequence, $x$ plays the role of $\xi$ in Section~\ref{sec:generalProgram}. 
	
	\paragraph{Numerical full measurement setup.} The full measurement setup considers uniformly weighted measurements taken at all inner vertices, \gls{ie}  
	\begin{equation}\label{eq:allLocations}
		{\Xi=\{x_n\}_{n=1}^{(N_x+1)^2}\backslash\partial X  \qquad \text{ with }\qquad |\Xi|=N=(N_x-1)^2. }
	\end{equation}

	\paragraph{Computation of $G_{x,:}$.} 
	For simplicity, we choose the background linearization parameter as the ground-truth $p_\ast = p_0$, and use both names interchangeably below. Evaluation of $\tilde \pi$ requires computation of the gradient $G_{x,:}= (\grd_pu_{p_\ast}(x))^\top$ for all $x\in\Xi$. In~\Cref{app:DerivationA(x)} we spell out the details of this computation using an adjoint method. For example, the $k$-th entry of the gradient reveals
	\begin{equation}\label{eq:gradient}
		[G_{x,:}]_{k} = [(\grd_p u_{p_\ast}(x))^\top]_{k} = \langle g^{(x)} ,  \phi_k u_{p_\ast}\rangle_{L^2(X)},
	\end{equation}
	where $g^{(x)} $ satisfies the adjoint equation
	\begin{align}\label{adjoint}
		-\Delta g^{(x)} + p_\ast g^{(x)}  = - \delta_{x} \text{ on } X, \quad g^{(x)} = 0\text{ on } \partial X.
	\end{align}

	\subsubsection{Optimal Sampling Distributions} \label{Sec:NumericalInvestigation}
	
	Usually the optimal sampling strategy $\tilde \pi$ is not available. As a numerical study, we nevertheless compute it by brute force and plot the result to examine its pattern and dependence on the ground-truth $p_\ast= p_0$ and correlated noise model $\Gamma$.

\begin{figure}
		\includegraphics[width = 0.24\textwidth]{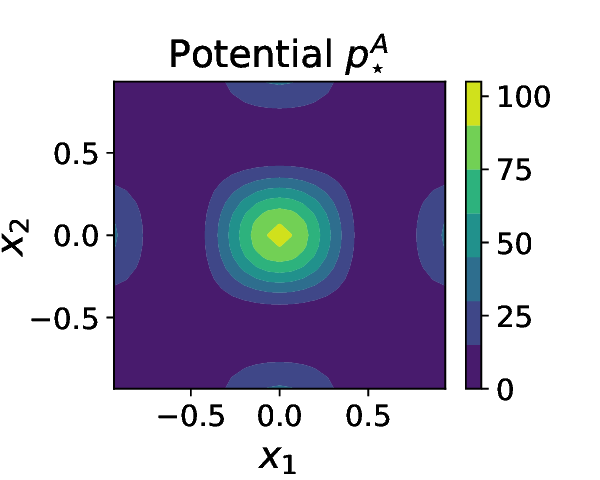}
		\includegraphics[width = 0.24\textwidth]{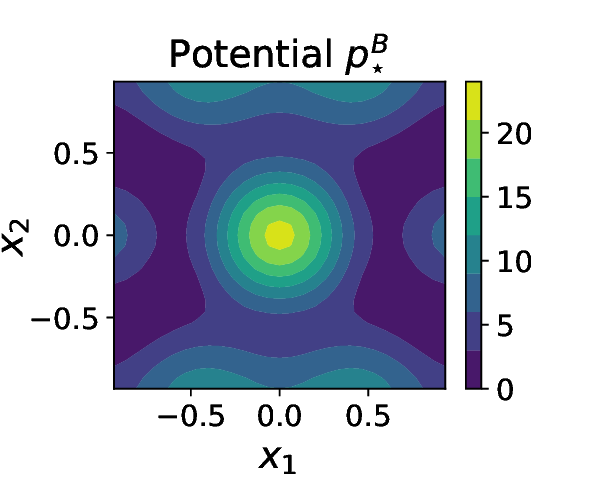}
		\includegraphics[width = 0.24\textwidth]{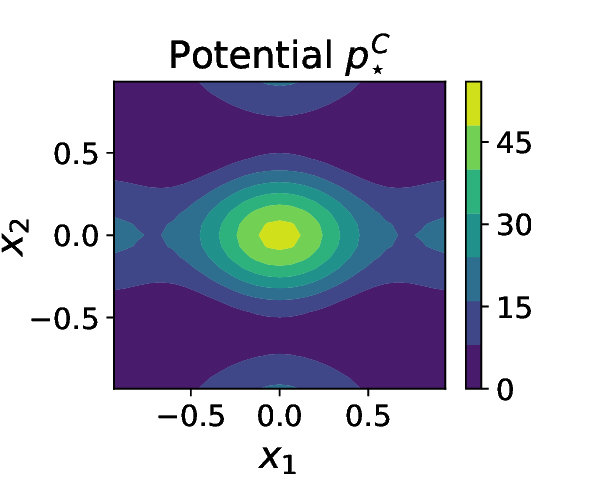}
		\includegraphics[width = 0.24\textwidth]{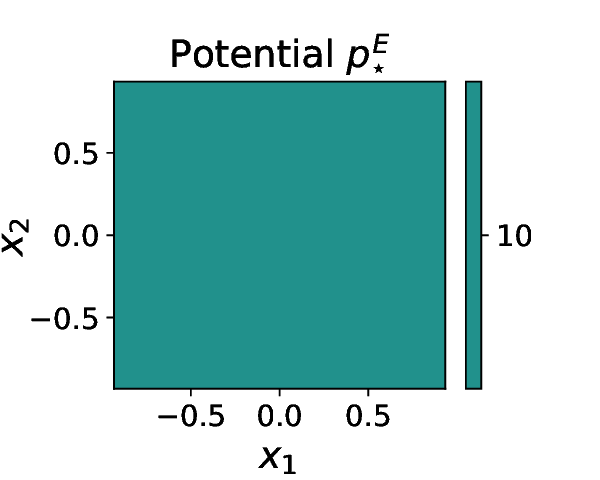}\\
		\includegraphics[width = 0.24\textwidth]{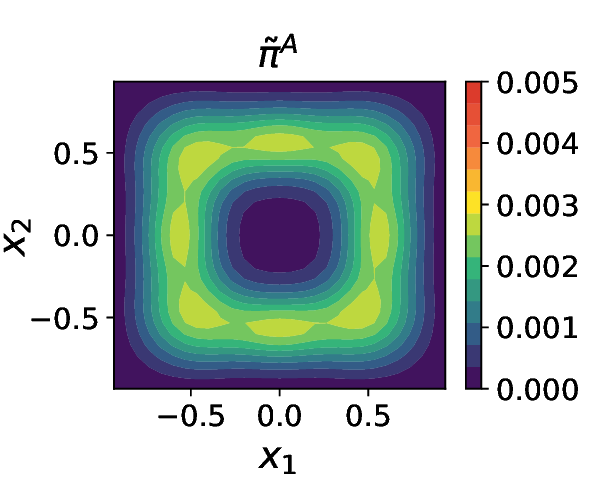}
		\includegraphics[width = 0.24\textwidth]{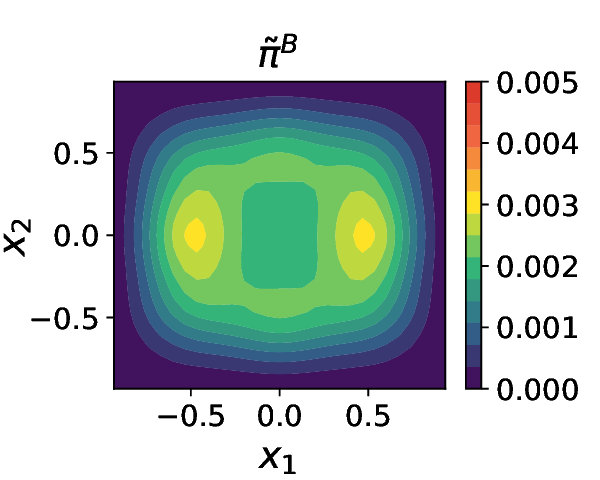}
		\includegraphics[width = 0.24\textwidth]{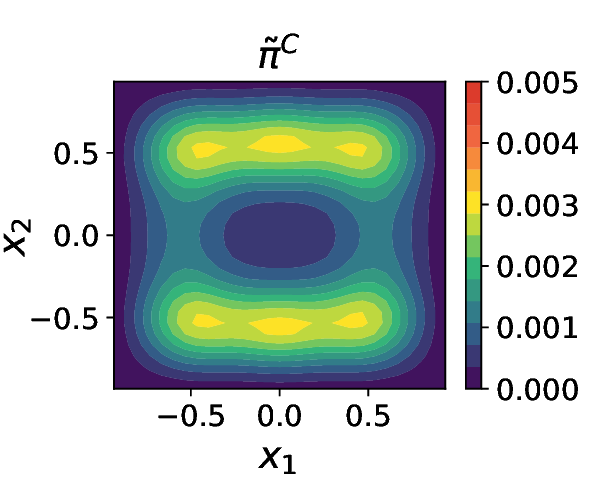}
		\includegraphics[width = 0.24\textwidth]{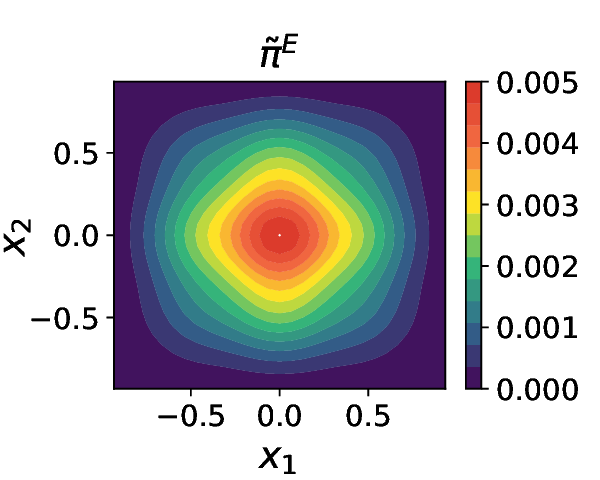}
		\caption{Top row shows four different  ground-truth media $p_\ast = p_0$, and the bottom row shows the optimal sampling distribution $\tilde\pi$ for  each of them.}
		\label{fig:piLandscapes2}
	\end{figure}

	\paragraph{ {Dependence on ground-truth $p_\ast$.}} In the four examples shown, we set the ground-truth parameters as in~\Cref{tab:groundtruthParams}, and  {uncorrelated noise $\eta(x) \sim \mathcal N(0,1)$ for all $x$}. The optimal sampling distribution in this case, according to~\eqref{eq:uncorrelatedSketching}, is  {$\tilde\pi(x)\propto \|G_{x,:}\|_2^2$, and this distribution is plotted in~\Cref{fig:piLandscapes2}.} The pattern of the optimal distribution shows significant dependence on the underlying ground-truth parameter $p_\ast$.
	
	\medskip 
	
	\begin{table}[H]
		\centering
		\tablinesep=2ex\tabcolsep=2pt
		\begin{tabular}{c||c|c|c|c}
			System $X$ & $A$ 
			&
			$B$ 
			& 
			$C$ 
			&  
			$D$ 
			\\\hline
			$p_\ast^X$&
			$
			\begin{pmatrix}
				13.6 & 10 &10\\
				10&10&10\\
				10&10&10
			\end{pmatrix}$
			& $
			\begin{pmatrix}
				5.856& 0.103& 3.168\\
				3.7441& 2.493& 1.124\\
				0.9902& 3.803& 0.846
			\end{pmatrix}$
			&  $
			\begin{pmatrix}
				11 &  8.889&  7.778\\
				6.667&  5.556&  4.444\\
				3.333& 2.222&  1.111
			\end{pmatrix}$
			& $
			\begin{pmatrix}
				10 & 0&  0\\
				0&  0&  0\\
				0& 0&  0
			\end{pmatrix}$
		\end{tabular}
		\caption{Test scenarios to study the optimal sampling strategy $\tilde \pi$. The $(i,j)$ entry of the matrix is the coefficient for $p_k$ with $(k_1 =i, k_2=j)$.}
		\label{tab:groundtruthParams}
	\end{table}

\begin{figure}[htb]
		\begin{adjustwidth}{-1cm}{-1cm}
			\centering
			\begin{subfigure}[c]{1.2\textwidth}
				\centering
				\includegraphics[width = 0.32 \textwidth]{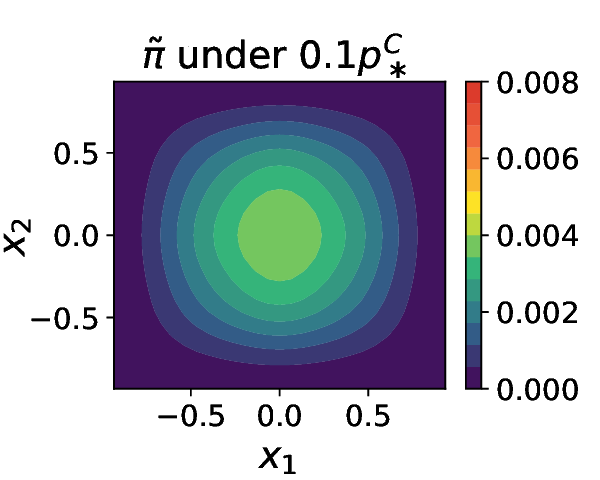}
				~ \includegraphics[width = 0.32 \textwidth]{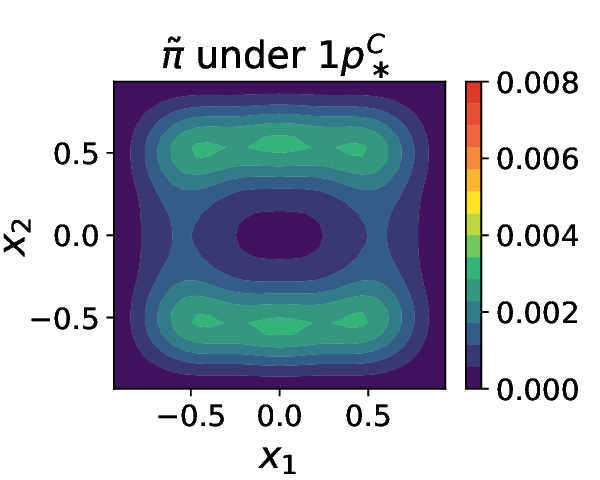}
				~ \includegraphics[width = 0.32 \textwidth]{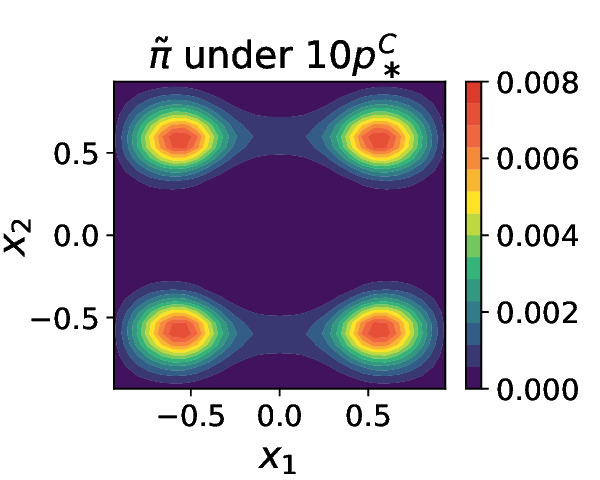}
				\subcaption{System C.}
			\end{subfigure}
			\begin{subfigure}[c]{1.2\textwidth}
				\centering
				\includegraphics[width = 0.32 \textwidth]{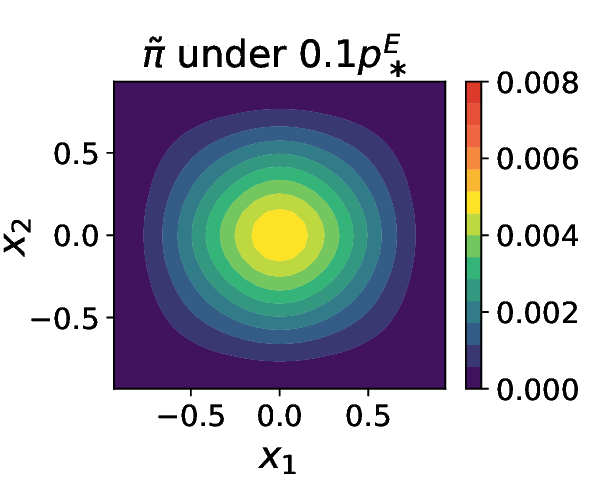}
				~ \includegraphics[width = 0.32 \textwidth]{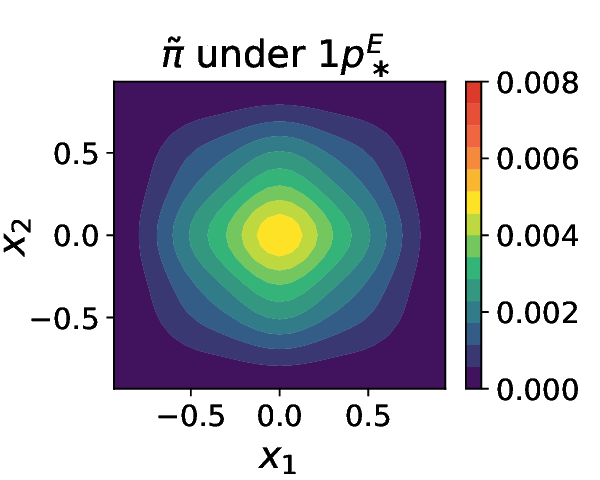}
				~ \includegraphics[width = 0.32 \textwidth]{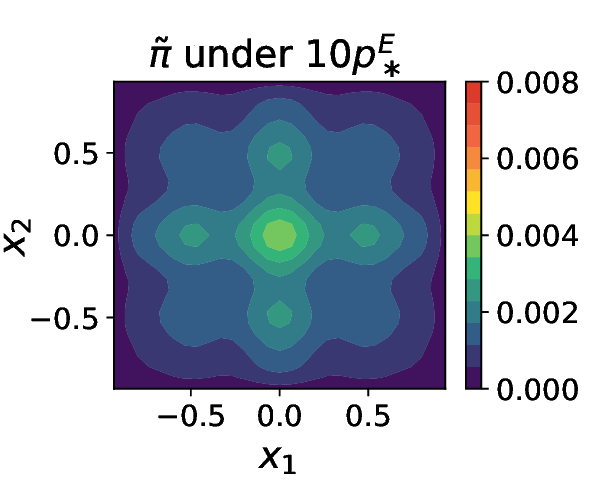}
				\subcaption{System D.}
			\end{subfigure}
			\caption{ Optimal importance sampling distributions $\tilde \pi$  for scaling parameters  $\alpha p_\ast$ with $\alpha= 0.1$ (left), $\alpha = 1$ (center) and $\alpha = 10$ (right). The ground-truth parameters $p_\ast$ from System C  and D from \Cref{tab:groundtruthParams} are taken.}
			\label{fig:piLandscapes:scales}
		\end{adjustwidth}
	\end{figure}

	{To investigate scaling effects,} we then scale the  { ground-truth} parameters by multiplying $p_\ast$ with a scaling parameter $\alpha {\in \{0.1, 1,10\}}$. We observe very different pattern{s} for  {the sampling distributions} $\tilde\pi$ when $\alpha$ is varied, as shown in~\Cref{fig:piLandscapes:scales}. The plots suggest that $\tilde\pi$ is centered in the middle when $p_\ast$ takes on small values, but develops interesting patterns when $p_\ast$ has a large scal{e}.
	
	{
		\paragraph{Dependence on noise model $\Gamma$.} In this experiment, we fix the ground-truth parameter to be $p_\ast^C$ as in~\Cref{tab:groundtruthParams}, and prescribe different correlated noise models $\eta \sim \mathcal N(0,\Gamma)$. The three noise models are:
		\begin{enumerate}[label=(\alph*)]
			\item Globally uniform strong interaction: $\Gamma^{-1}_{i,j} = -\delta_{i\neq j} + (N+1) \delta_{i=j}$
			\item Local medium interaction: only the four neighboring measurement locations interact with the current location and determine it partially: $$\Gamma^{-1}[i,j] = -\delta_{x_j\in \overline{B_{dx}(x_i)}} +  (4+1) \delta_{i=j}$$
			\item Midrange spatial interactions: interaction is stronger for closer neighbors and close to the boundary of the domain 
			\begin{align*}
				&\Gamma^{-1}[i,j] = -e^{-(x_j-x_i)^2/(2\gamma^2)}(x_j^2+0.5)(x_i^2+0.5), \quad \text{ and}\\
				&\Gamma^{-1}[i,i] = -\sum_{j\neq i}\Gamma^{-1}[i,j] +1.
			\end{align*}
	\end{enumerate}}
	
\begin{figure}
		\centering
		\includegraphics[width=0.5\linewidth]{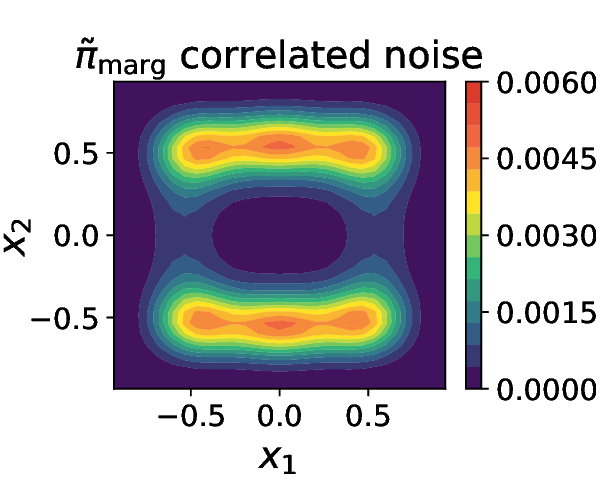}~\includegraphics[width=0.5\linewidth]{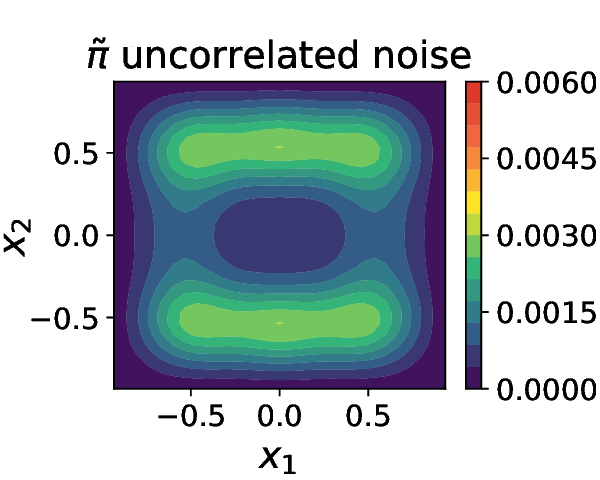}\\
		\includegraphics[width=0.5\linewidth]{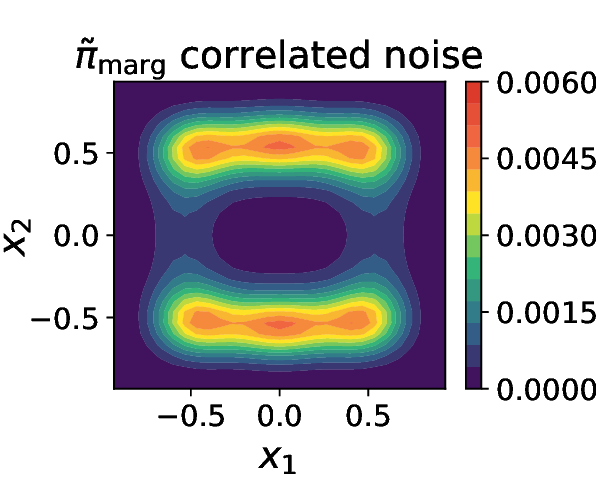}~\includegraphics[width=0.5\linewidth]{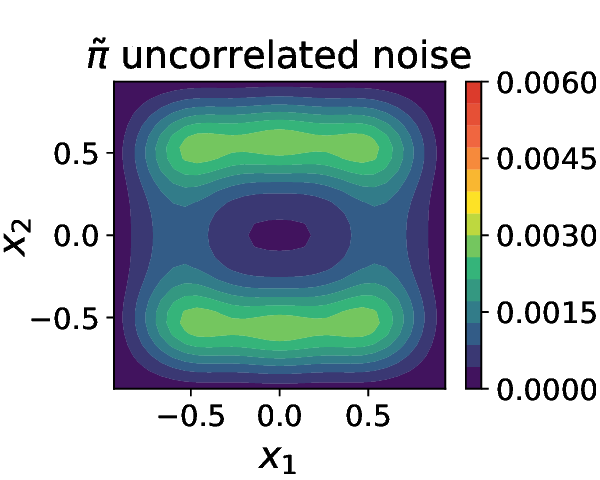}\\
		\includegraphics[width=0.5\linewidth]{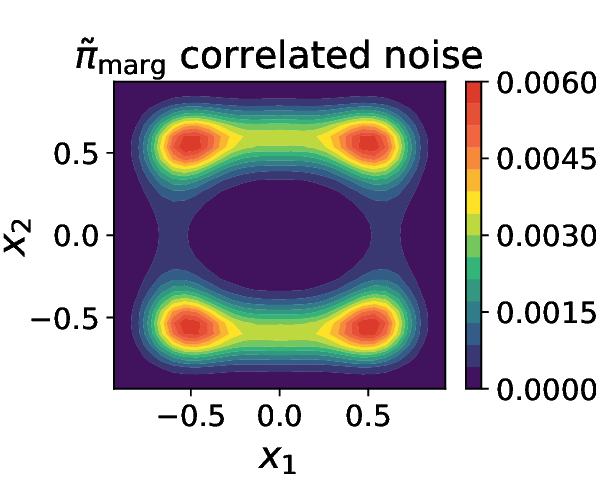}~\includegraphics[width=0.5\linewidth]{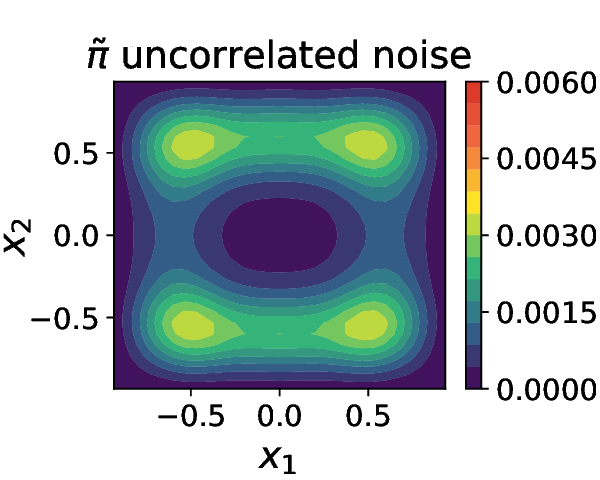}
		\caption{Sampling distributions for System C under different noise correlations (left column): global interaction (first row), local interaction (second row) and midrange spatially aware  interactions (third row).  They are compared to  the respective uncorrelated noise model versions that attain the same variance (right column).}
		\label{fig:CorrelatedNoise}
	\end{figure}

	{With different choices of $\Gamma$, we plot in \Cref{fig:CorrelatedNoise} the resulting optimal distribution marginalized on its first variable: $\tilde \pi_{\mathrm{marg}}(x) = \sum_y \tilde \pi(x,y)$. Different correlation and noise structures give rise to very different optimal sampling distributions $\tilde \pi$. Their corresponding uncorrelated versions are obtained by removing all off-diagonal entries in $\Gamma$.}

	\subsubsection{Effect of Sampling on sensor locations and \gls{fim}} \label{Sec:NumericalEvidence}

	As a proof of concept, we now study the performance of our sampling algorithms for  {the} recovery of  {suitable} sensor locations. Both the inverse condition number $c_{\mathrm{inv}}$ and the minimal eigenvalue $\lambda_{\min}$ of the {\gls{fim}} will be tracked and reported.  {We regard our final sensor selection as a success if the subsampled \gls{fim} has $c_{\mathrm{inv}}\approx 1$ and $\lambda_{\min}\gg0$.} 
	
	{Experiments in this section assume the ground-truth parameter is $p_\ast^C$ from \Cref{tab:groundtruthParams}, and the data is generated accordingly.}
	
\begin{figure}
		\centering
		\includegraphics[width = 7cm]{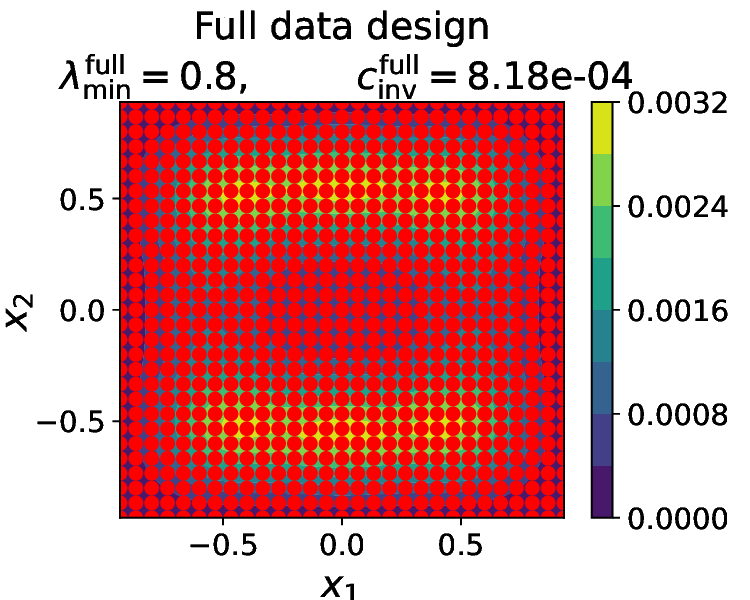}
		\caption{Full-Data Setup: Measurement locations (red dots) are located in all grid points. The optimal importance sampling distribution $\tilde\pi$ is drawn in the background.}
		\label{fig:Sampling_FullData}
	\end{figure}

	{The full dataset is given in terms of measurements everywhere on a grid with $N_x=30$ cells in every direction, resulting in {$N=(N_x-1)^2=841$} sensor locations. These are depicted as red dots in~\Cref{fig:Sampling_FullData}. The background shows the optimal distribution $\tilde \pi$. In this full-data setting, the \gls{fim} attains an inverse condition number   of $c_\inv^\full = 8e{-4}$ and a minimum eigenvalue of $\lambda_{\min}^{\full} = 0.8>0$, and the problem is locally {sensitive}.}
	
\begin{figure}
		\centering
		\includegraphics[width=0.5\linewidth]{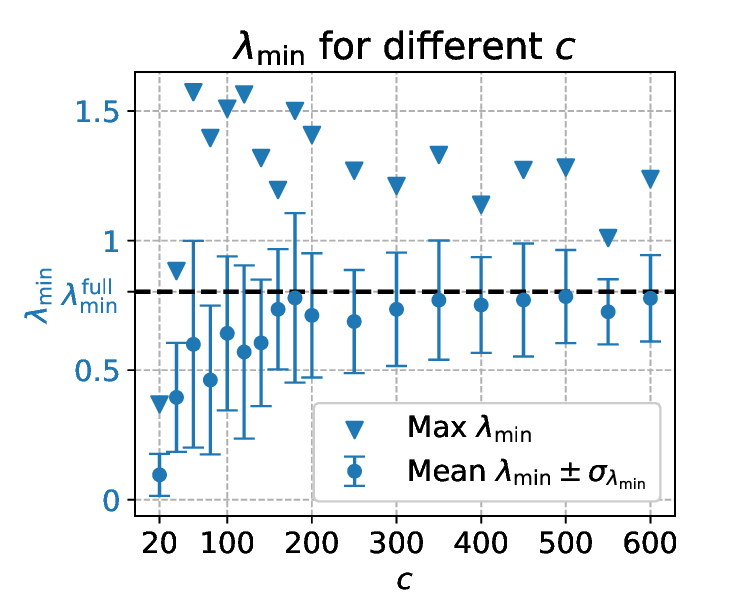}~\includegraphics[width=0.5\linewidth]{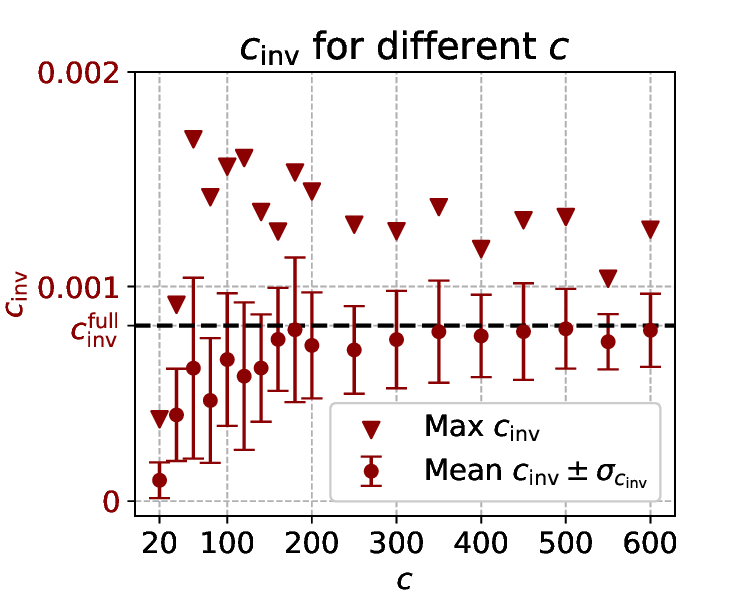}
		\caption{Mean, standard deviation range and maximum observed minimal eigenvalue (left) and inverse condition number (right) of the \gls{fim} under $40$ independent realizations of $c$ samples from $\tilde \pi$.}
		\label{fig:SamplingStatistics}
	\end{figure}

	\paragraph{Effect of Importance Sampling.}
	{In \Cref{fig:SamplingStatistics}, we plot statistical behaviour of the minimal eigenvalue and inverse condition number for the subsampled \gls{fim} using $c$ samples drawn from $\tilde \pi$. For each sample size $c$, we perform  40 independent sampling realizations and plot the  mean (dot), the maximum value (triangle) and the $\pm1$ standard deviation range around the mean (whiskers).
	}
	
	{As predicted by theory, the mean stabilizes around the value for the full-data setting, with shrinking variance as $c$ grows. This reduced variability is also reflected in a decline of the maximum values. 
	}
	
	\paragraph{Effect of Pipeline \Cref{alg:SampleDesigns_EnsembleGreedy}.}
	{Throughout this part, we employ}   {the ensemble-based sampling strategies as described in \Cref{ssec:PracticalConsiderations},  with early stopping based on the condition number. }

\begin{figure}
		\centering
		\includegraphics[width = 0.32\textwidth]{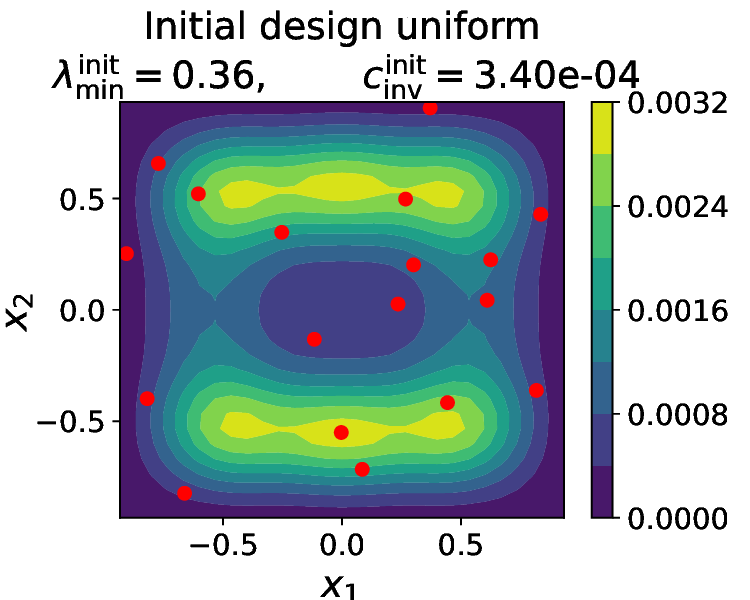}~
		\includegraphics[width = 0.32\textwidth]{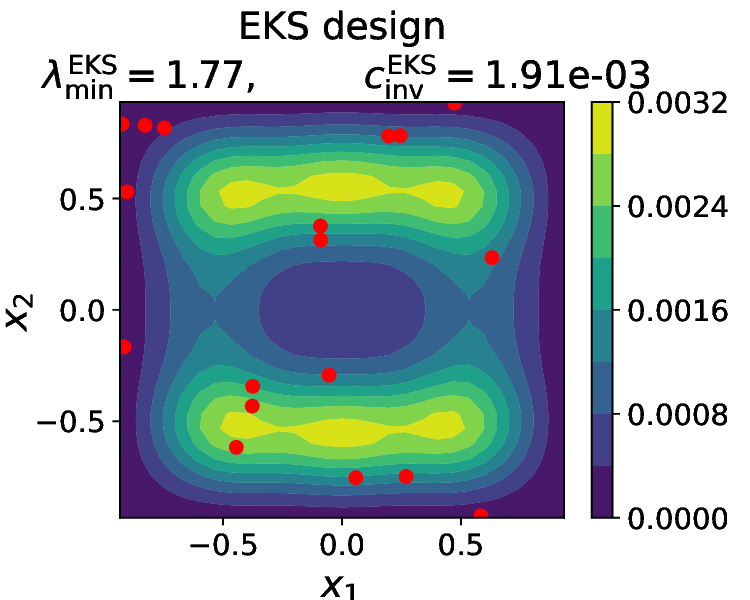}~
		\includegraphics[width = 0.32\textwidth]{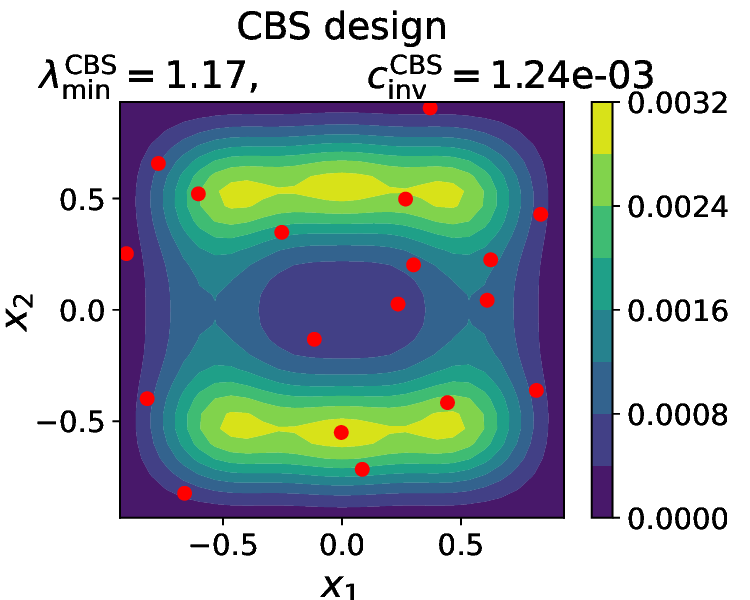}
		\caption{Uniformly distributed initial guess (left) of the distribution of the sensors (red dots) in the domain $X$, where the optimal importance sampling distribution is drawn in the background. Application of  {early-stopped} \gls{eks} (middle) and \gls{cbs} ({right})  {sampling from the sensitivity distribution $\tilde \pi$}
			changes the sensor distribution and dramatically increases the condition number  and minimum eigenvalue of the respective \gls{fim}.} 
		\label{fig:Sampling_UniformData}
	\end{figure}

	{For} down-sampling, we allow only $c=18 = 2K$ sensor locations. {A naive initial guess could be to place  sensors {uniformly} over $\Xi$}  { with equally weighted data points. 
		The resulting sensor configuration is plotted in the left panel of \Cref{fig:Sampling_UniformData}. Interestingly, this configuration  already performs reasonably well, attaining a minimum eigenvalue and conditioning $(\lambda_{\min}^{\init,u},c_{\inv}^{\init,u})$ of $(0.36,3.4e{-4})$ close to the full-data setting.} Note that the {optimal} sampling distribution $\tilde \pi$ is bounded from above by $0.0031$. Hence, the uniform distribution $\frac{1}{N}\approx0.0012$ is, in this particular case, a good approximation (with $\beta \leq 0.383$), resonating our discussion in Section~\ref{ssec:earlystopping}.
	{We then invoke our \Cref{alg:EarlyStoppingSampling} for importance sampling and reweighting with $\tilde \pi$ through \gls{eks} or \gls{cbs} with  {early stopping}, for $25$  {iterations}. Both} move the samples to new locations, and further improve the conditioning of the  \gls{fim}: $(\lambda_{\min}^\EKS,c_\inv^\EKS)$ and $(\lambda_{\min}^\CBS,c_\inv^\CBS)$ become $(1.77,1.91e{-3})$ and $(1.17,1.24e{-3})$ respectively,  {as summarized in \Cref{tab:DesignSamplingQ9}. }

	\begin{table}
		\centering
		\begin{tabular}{lcc}
			Design $D$&$\lambda_{\min}^D$&$c_{\mathrm{inv}}^D$ \\\hline\hline
			full data $D^{\mathrm{full}}$ &$ 0.8$& $8.18\cdot 10^{-4}$\\\hline
			uniform initial guess 
			&$0.36$ & $3.4\cdot 10^{-4}$ \\
			\gls{eks} sample 
			&  $1.77$&$1.91\cdot 10^{-3}$\\
			\gls{cbs} sample
			& $1.17$&$1.24\cdot 10^{-3}$
			\\\hline
			normal initial guess
			&$1.48\cdot 10^{-4}$ &  $1.54\cdot 10^{-7}$\\
			\gls{eks} sample 
			&  $2.06$&$2.25\cdot 10^{-3}$\\
			\gls{cbs} sample 
			& $1.41$& $1.56\cdot 10^{-3}$\\  
		\end{tabular}
		\caption{Comparison of {sensitivity} measures  associated with different designs. Rows below an initial guess  refer to sampling starting from this  initial design.}
		\label{tab:DesignSamplingQ9}
	\end{table}

\begin{figure}
		\centering
		\includegraphics[width = 0.32\textwidth]{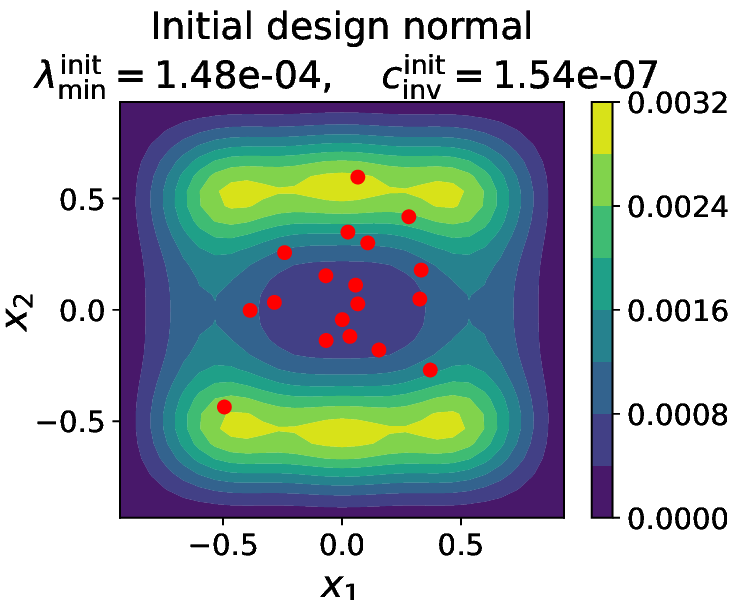}~
		\includegraphics[width = 0.32\textwidth]{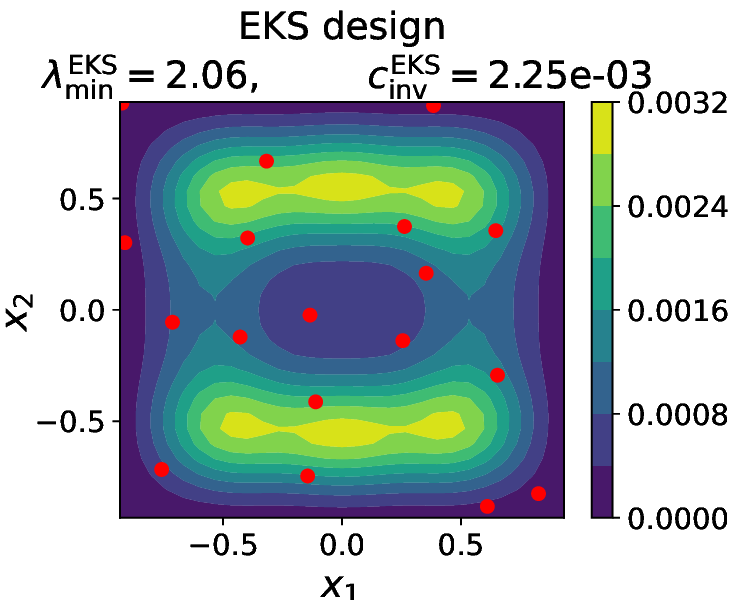}~
		\includegraphics[width = 0.32\textwidth]{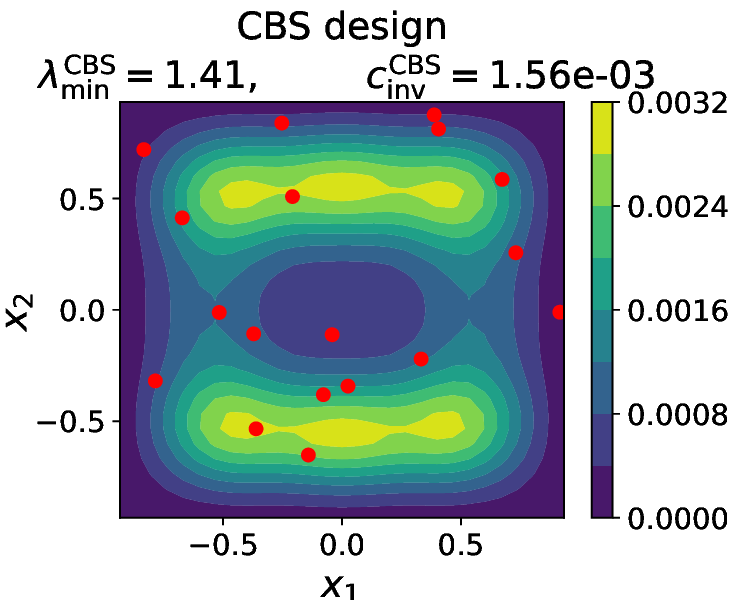}
		\caption{Red markers demonstrate the location of the sensors, with the background plotted as the optimal distribution. The left panel shows the distribution of the initial samples. The middle and the right panel show, respectively, the \gls{eks} and \gls{cbs} samples after 25 iterations. The minimum eigenvalues of the \gls{fim} change from $1.48e{-4}$ to $2.06$ and $1.41$ and the inverse condition numbers from order $1e{-7}$ to $1e{-3}$,  ensuring local {sensitivity}.}
		\label{fig:Sampling_NormalData}
	\end{figure}

\begin{figure}
		\centering
		\includegraphics[width = 7cm]{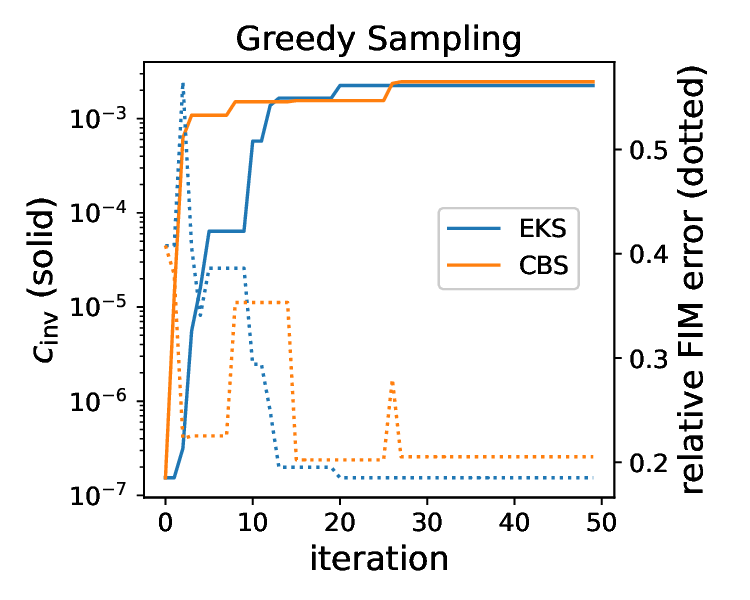}
		\caption{
			Evolution of {the inverse condition number} (solid lines) and {relative Frobenius norm error}  of the down-sampled \glspl{fim} from the full-data setup (dashed lines)  {under }  \gls{eks} (blue)  {and} \gls{cbs} (orange)  {sampling} with early stopping.   {The same centered  Gaussian $\mathcal N(0,\frac 1 4 \operatorname{Id}_2)$ distributed initial sensor locations} are shared across the sampling methods.} \label{fig:evolveLambdaminEvolution_NormalData}
	\end{figure}

	{Our method unveils its full capability  when we consider  a poorly chosen initial guess for the sensor locations, given by a narrowly centered normal distribution {$\mathcal N(0, \frac{1}{4}\mathrm{Id}_2)$. }
		The resulting sensor configuration is plotted in the left panel of \Cref{fig:Sampling_NormalData}, and attains a substantially worse reconstruction setting}, achieving $c_{\inv}^\init = 1.54e{-7}$ and $\lambda_{\min}^{\init} = 1.48e{-4}$. We once again apply~\Cref{alg:EarlyStoppingSampling} that integrates \gls{eks} or \gls{cbs}. The conditioning of the weighted \gls{fim} is significantly improved, with $(c_\inv^\EKS,\lambda_{\min}^\EKS) = (2.25e{-3}, 2.06)$, and  $(c_\inv^\CBS,\lambda_{\min}^\CBS) = (1.56e{-3},1.41)$, respectively.  The approximation errors in the \gls{fim} along  {the ensemble sampling iterations}  are plotted in~\Cref{fig:evolveLambdaminEvolution_NormalData}.

	We observe that  {in both cases} the inverse condition number  and  the minimum eigenvalue  of the \gls{fim} corresponding to designs generated by \gls{eks} and \gls{cbs}  {exceed}  those obtained using the full dataset. {This behaviour is targeted by  the pipeline through the early stopping at favourable configurations and importance-reweighting, compare \Cref{fig:SamplingStatistics}.} 
	This mechanism focuses our experimental design strategy on selecting informative data points and emphasizing them more heavily, thereby amplifying the overall information content.


	\subsubsection{Effect of Sampling on the square loss function}

	The sensitivity of the data \gls{wrt} the parameter can also be observed by the strong convexity of the quadratic cost function $\C(p) = \| y(\cdot)-\F(\cdot,p)\|^2_{2}$. In what follows, we visualize the landscape of this cost function across the parameter space for different experimental designs, in order to assess the impact of our sampling strategy on data sensitivity and the difficulty of the full{y} nonlinear inversion problem.
	
\begin{figure}
		\begin{adjustwidth}{-1cm}{-1cm}
			\centering
			\includegraphics[width = 6cm]{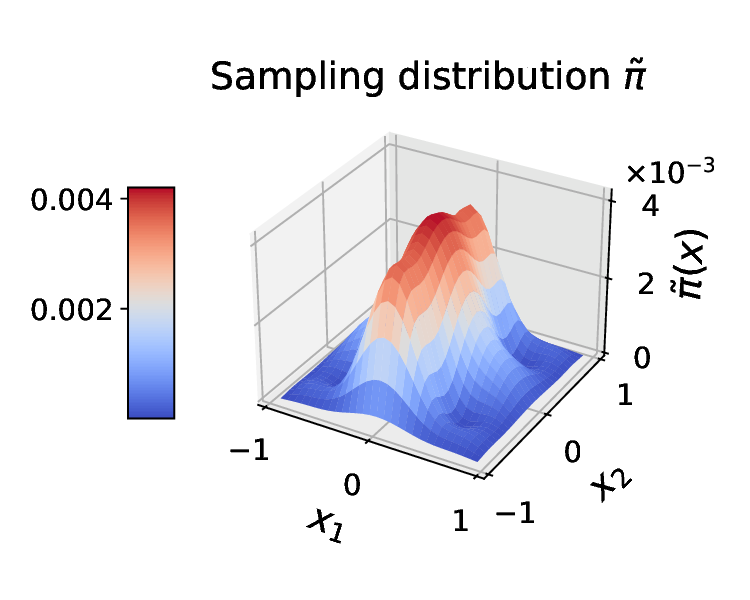} ~\includegraphics[width = 6cm]{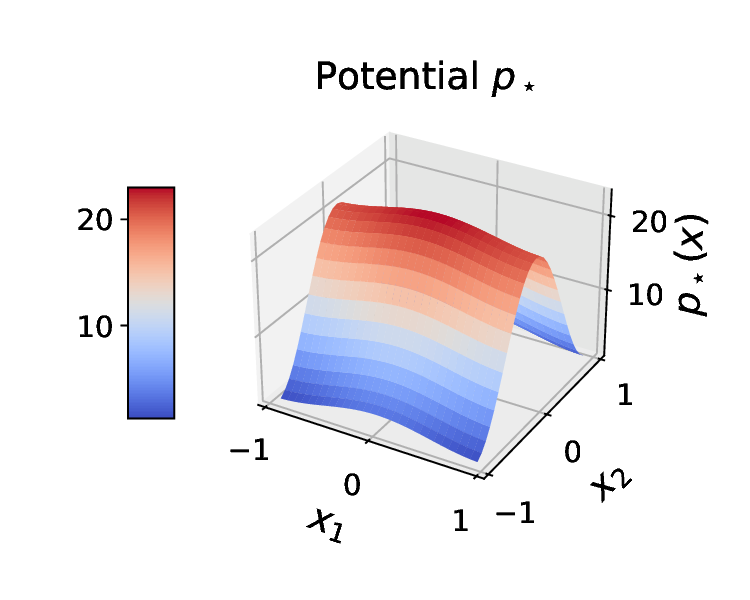}
		\end{adjustwidth}
		\caption{Optimal importance sampling distribution $\tilde\pi$ (left) and shape of the ground-truth parameter $p_\ast$ (right) in the two-dimensional  {parameter}  setting. }
		\label{fig:LossLandscape:Parameter}
	\end{figure}

	For visualization, we  confine ourselves to a  {(K=2)}-dimensional admissible set with $\mathcal{A} = \{p:X\to {\rr_0^+}\mid p(x) = p_1\cos(x_1)+ p_2\cos(x_2)+12\}$ and  {set} the ground-truth parameter  {to} $p_0(x_1,x_2)= 1\cos(x_1)+ 10\cos(x_2)+12$, set $p_\ast = p_0$ and work with noise free synthetic data. The profile of $p_\ast$ and the optimal importance sampling distribution are depicted in \Cref{fig:LossLandscape:Parameter}. The scaling for $p_\ast$ in the $x_1$ and $x_2$ direction is very different, and $p_\ast$ changes its profile in $x_2$ direction significantly more. This is in alignment with the extension of the sampling probability.
	
\begin{figure}
		\centering
		\includegraphics[width =  0.45\textwidth]{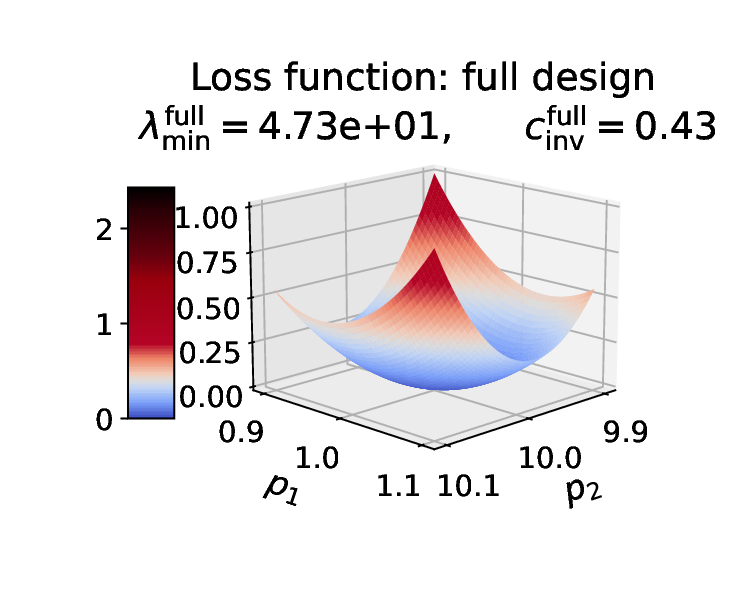} 
		~  \includegraphics[width =  0.45\textwidth]{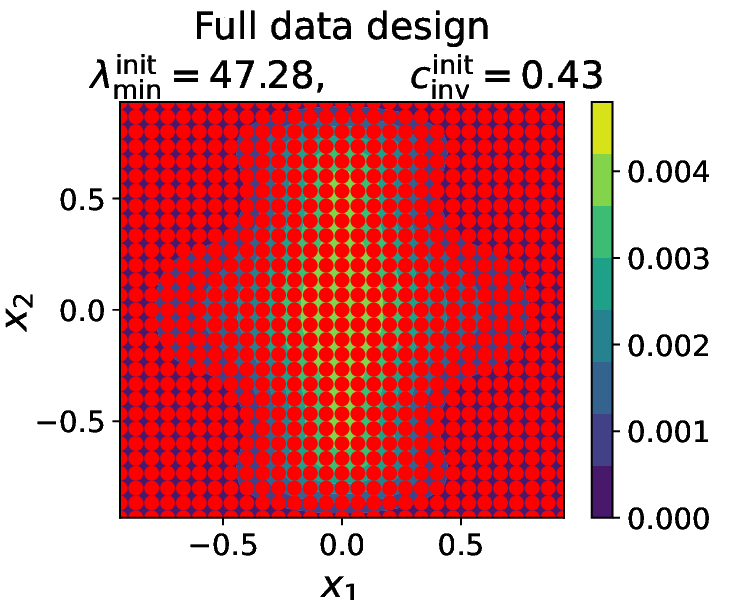}\\ 
		\includegraphics[width =  0.45\textwidth]{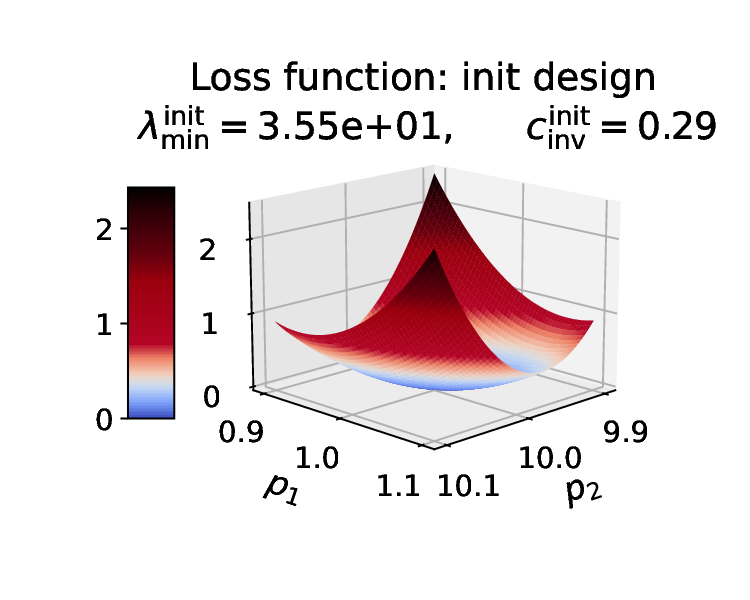} 
		~  \includegraphics[width =  0.45\textwidth]{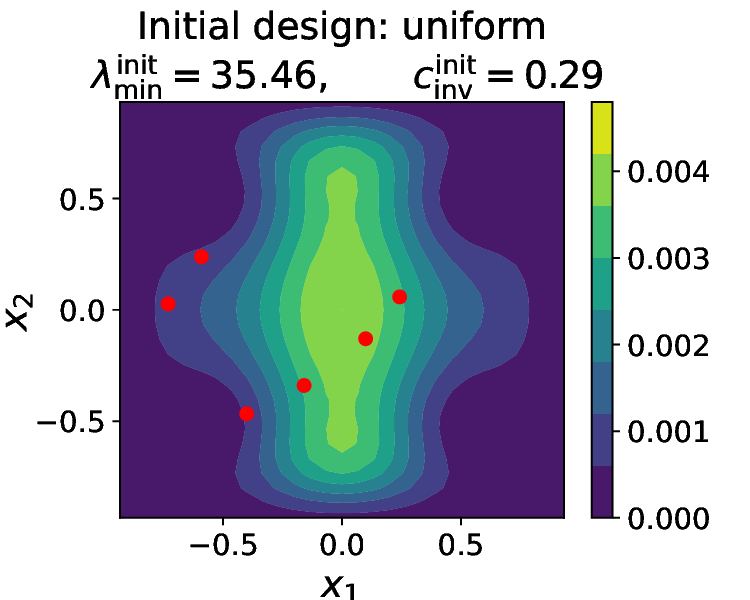}
		
		\caption{Loss landscapes (left) for different sensor locations (right): full-data setup (first row) and  {uniformly} distributed initial sensor locations (second row).}
		\label{fig:LossLandscapes_fullInit}
	\end{figure}

\begin{figure}
		\includegraphics[width = 0.45\textwidth]{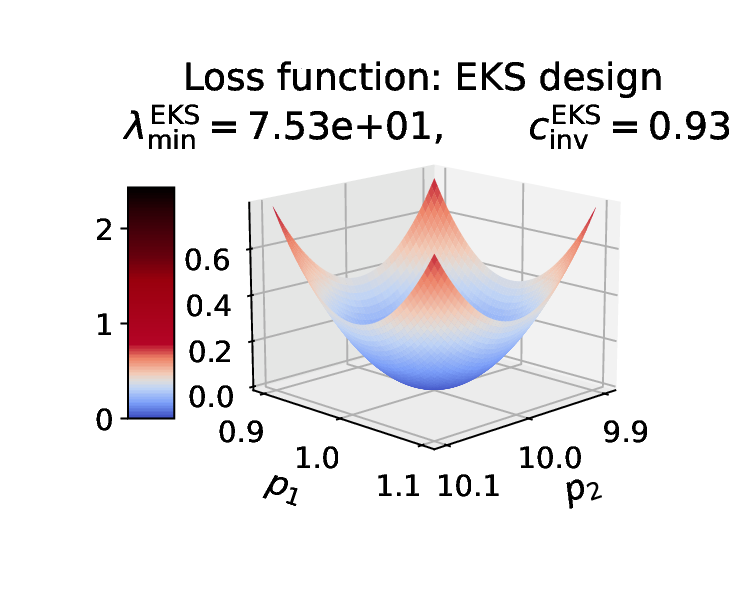}
		\includegraphics[width = 0.45\textwidth]{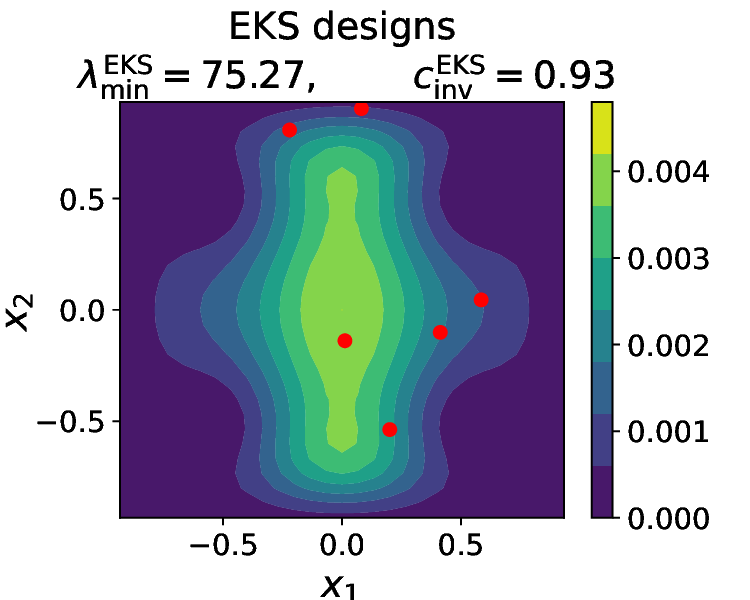}\hspace{1cm}\\
		\includegraphics[width = 0.45\textwidth]{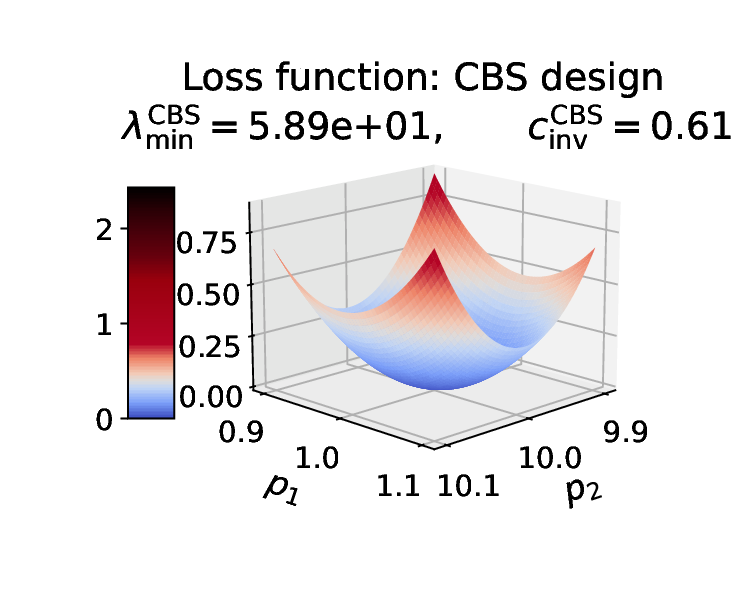}
		\includegraphics[width = 0.45\textwidth]{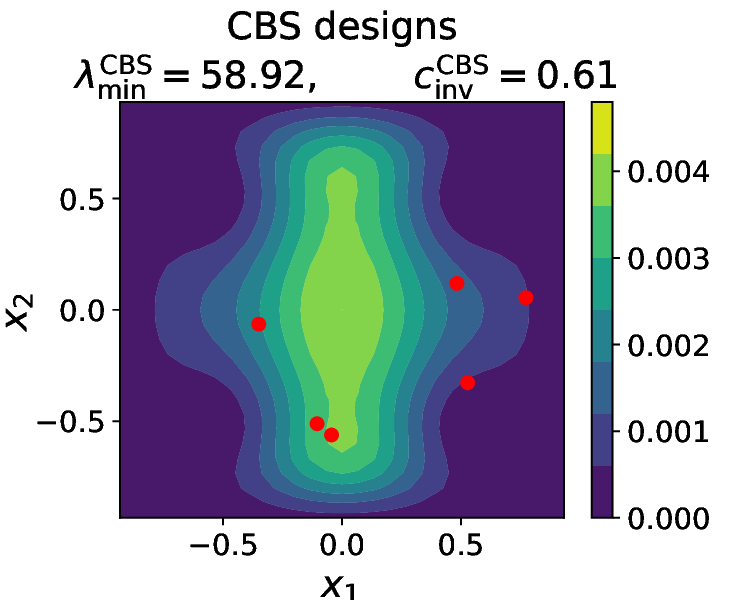}
		\caption{Quadratic loss landscapes (left)  {corresponding to different sensor  locations (right) sampled from the sensitivity-based sampling distribution $\tilde \pi$ (background) with \gls{eks} (upper row) or  \gls{cbs} (lower row) under early stopping.}}
		\label{fig:LossLandscapes_Sampling}
	\end{figure}

\begin{figure}
		\centering
		\includegraphics[width =  0.45\textwidth]{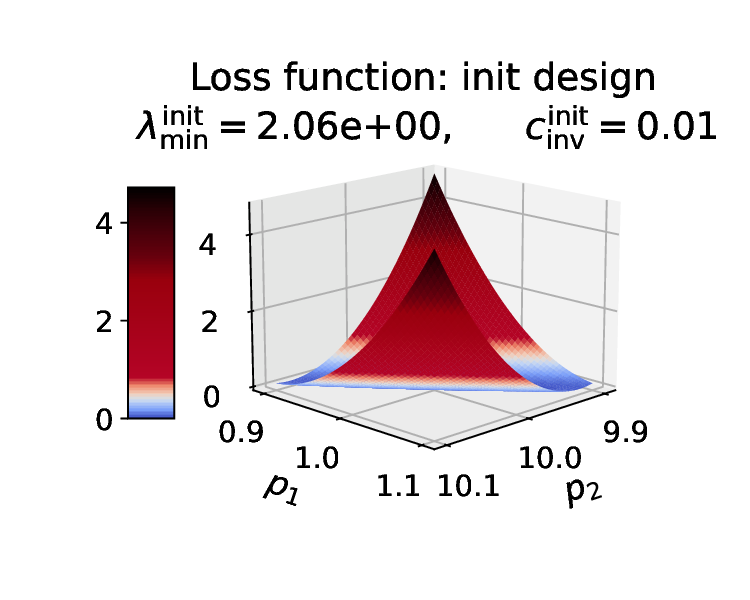} 
		~  \includegraphics[width = 0.45\textwidth]{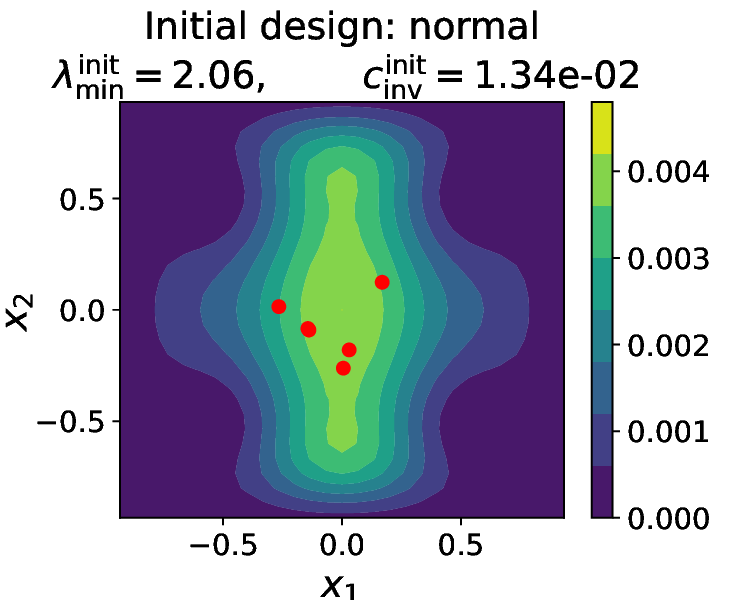}\\
		\includegraphics[width =  0.45\textwidth]{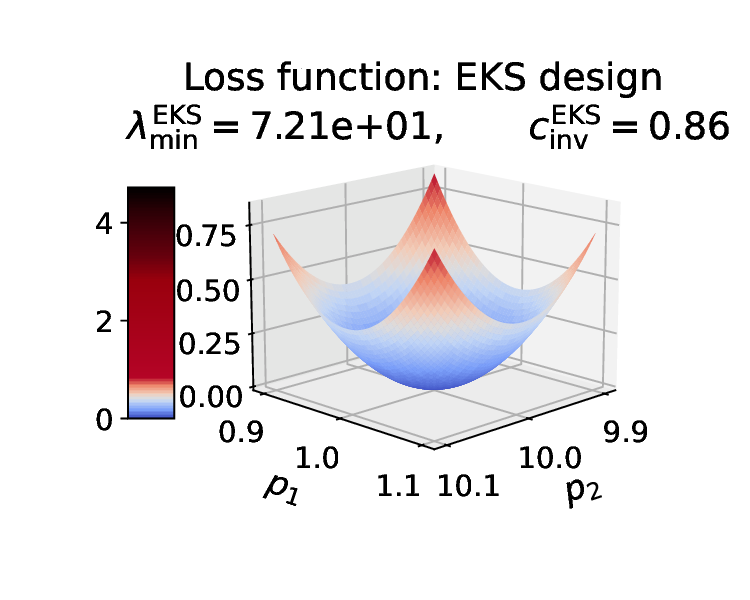}
		~\includegraphics[width =  0.45\textwidth]{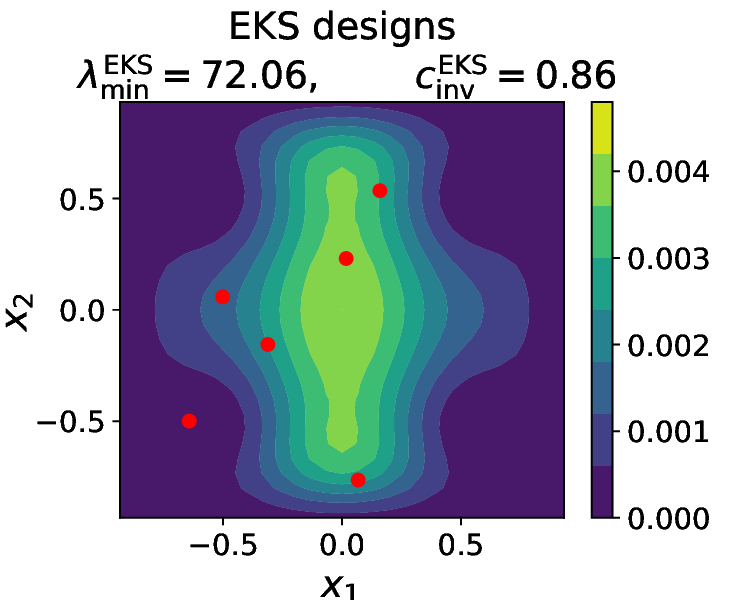}
		\caption{Loss landscapes (left) for different sensor locations (right):  {normally} distributed initial sensor locations (first row) and  sensitivity-based sensors sampled using EKS (second row).}
		\label{fig:LossLandscapes_Sampling_normal}
	\end{figure}

	When the full dataset is used, the loss function is convex, indicating the full dataset contains sufficient information for the recovery, with a  conditioning of $c_{\inv}^\full = 0.43$ and a  minimum eigenvalue being $\lambda_{\min}^\full = 47.3$, as shown in \Cref{fig:LossLandscapes_fullInit}.  {An initial choice of } $8$   { uniformly 
	} distributed sensor locations shows slightly reduced convexity in the landscape of the objective function, and the inverse condition number becomes $0.29$, with a minimum eigenvalue of $35.5$.  {\gls{eks} and \gls{cbs}} sampling with   {early stopping}  based on the condition number  enhances both the convexity and the conditioning,  {yielding landscapes whose convexity exceeds that of the full-data setting}, as Figure~\ref{fig:LossLandscapes_Sampling} shows.
	{Again, this improvement is more significant when it is compared with a poor initial sensor choice is poor, such as  } $8$   { normally $ \mathcal N(0, \frac{1}{4}\mathrm{Id}_2)$
	} distributed sensor locations as summarized  in   Figure~\ref{fig:LossLandscapes_Sampling_normal}.

	\subsection{Source Term Design}
	In our second set of experiments, we  { consider the source term as an additional design choice.} In particular, we  {prescribe its form as}:
	\[
	\gamma(x) = \gamma_1 x_1 + \gamma_2 x_2 + 10\,,\quad\text{with}\quad \vec \gamma = (\gamma_1, \gamma_2)\in [-2,2]^2\,.
	\]
	Similar to the previous example, the possible measurements are the solution evaluated at points $u^\gamma(x)$. The entire forward map is:
	\[
	\hat \F(x,\vec{\gamma}, p)= u^{\vec\gamma}(x)\,,\quad\text{and}\quad \hat y(x,\vec\gamma)= \hat \F(x,\vec\gamma, p_\ast)+ \eta(x, \vec\gamma)\,.
	\]
	The flexibility of $x$ and $\vec{\gamma}$ means we have a four-dimensional  design space: $\hat \Xi = \Xi \times  {\Omega}$  {for a fine discretization  $\Omega = \{-2+\frac{4i}{N_\gamma}\}_{i= 0,...,N_\gamma}^2$ of $[-2,2]^2$ with $N_\gamma \in \mathbb N$}. We  fix the parameter dimension to $K=9$ again.  {The fine discretization of the $\vec \gamma$ space $\Omega$} prevents us from computing the full landscape, as well as the normalization constant of $\tilde \pi$ exactly,  and we  use this  as an example to demonstrate our method in this setting.

	{At the initial stage, we place in total $18=2K$ combinations of source and sensor over the entire $\hat \Xi$ plane in a uniform manner. As in the case for a fixed source term, this initialization already yields a local sensitivity value of $2.25e{-4}$ for the inverse {\gls{fim}} conditioning. Running~\Cref{alg:EarlyStoppingSampling} provides slightly better conditioning, with both \gls{eks} and \gls{cbs} provide $c_{\inv}$ values exceeding $1e{-3}$. See \Cref{tab:DesignSamplingSource}.}

\begin{figure}[htb]
		\centering
		\includegraphics[width=0.32\linewidth]{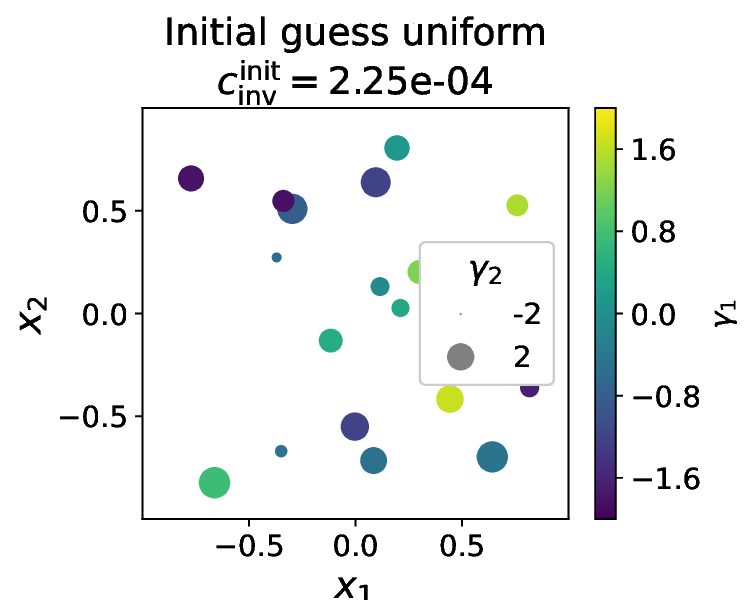}~
		\includegraphics[width=0.32\linewidth]{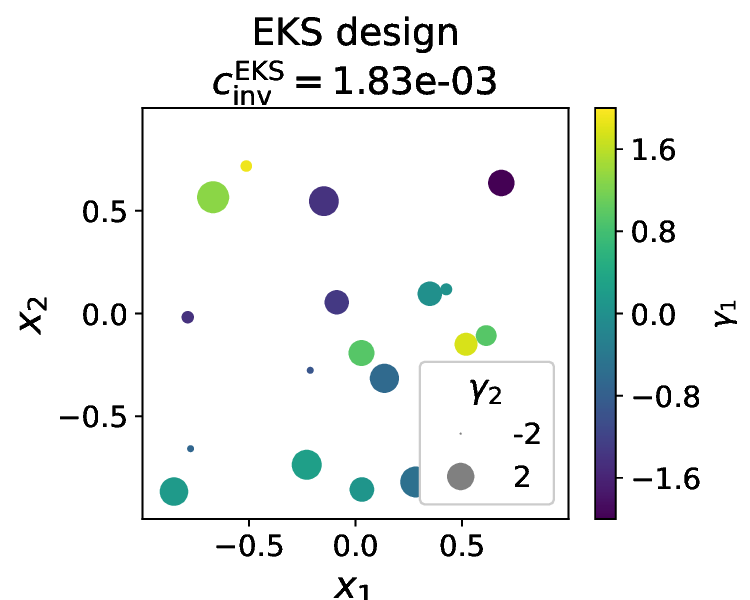}~
		\includegraphics[width=0.32\linewidth]{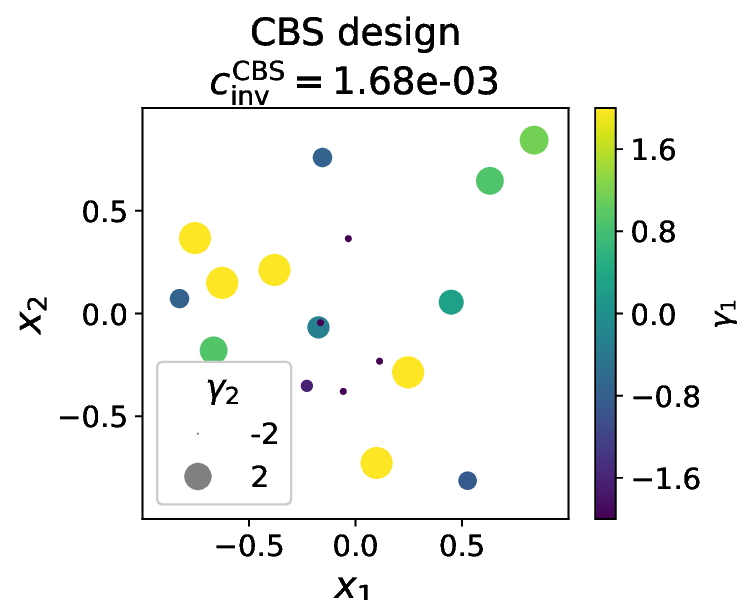}\\
		\caption{  {Three} different designs, characterized by their sensor locations given by the dot locations, and $\gamma_1,\gamma_2$ values  encoded in  colour and size of the dots, together with their {sensitivity} measures:   normally distributed initial sensor location guess with uniformly distributed $\gamma_1,\gamma_2$ (left),
			{early-stopped} \protect\gls{eks} (middle) and \protect\gls{cbs}  (right)  {sampling} \protect\gls{wrt} the rescaled sensitivity-based sampling distribution, after 60 iterations each.  }
		\label{fig:SamplingSource_Normal}
	\end{figure}

	{If, instead, sensor locations are initially chosen according to a very concentrated normal distribution {$\mathcal N(0,\frac 1 4 \operatorname{Id}_2)$}, then~\Cref{alg:EarlyStoppingSampling} returns better conditioning and pushes these samples to the wider domain, under  both  samplers \gls{eks} and \gls{cbs}, This is} illustrated in  \Cref{fig:SamplingSource_Normal}.  {The resulting sensitivity values are} summarized in \Cref{tab:DesignSamplingSource}.

	\begin{center}
		\begin{table}[H]
			\centering
			\begin{tabular}{lcc}
				Design $D$
				&$c_{\mathrm{inv}}^{D}$ \\\hline\hline
				uniform initial guess 
				& $2.25\cdot 10^{-4}$ \\
				\gls{eks} sample
				&$1.83\cdot 10^{-3}$\\
				\gls{cbs} sample 
				&$1.68\cdot 10^{-3}$ 
				\\\hline
				normal initial guess 
				&  $3.02\cdot 10^{-9}$\\
				\gls{eks} sample 
				&$1.68\cdot 10^{-3}$\\
				\gls{cbs} sample 
				& $8.75\cdot 10^{-4}$        
			\end{tabular}
			\caption{Comparison of the local {sensitivity}, measured through the inverse FIM conditioning, resulting from different design strategies. Rows below an initial guess refer to sampling starting from this  initial configuration.}
			\label{tab:DesignSamplingSource}
		\end{table}
	\end{center}
	
	\section{Discussion}\label{sec:Discussion}

	In this work, we examine qualitative experimental design of a parameterized inverse problem. By reformulating the preservation of  {data sensitivity} as the preservation of the positivity of the \gls{fim}, we examine the possibility and strategy for down-sampling the data. This problem is further reduced to matrix sketching, where a well-studied sketching algorithm from \gls{rnla} becomes handy. The sample size depends on a sampling distribution that reflects the structure of the forward problem, and ensemble-based sampling algorithms such as  \gls{eks} and \gls{cbs} are implemented.
	
	The general program described in this article can be applied to a variety of experimental design / data selection tasks arising from inverse problems. As a proof of concept, we provide a numerical test using  {the} Schr\"odinger equation as the forward model. The optimal distribution is problem-dependent and is typically unavailable. In various applications, knowledge of the forward model can be used to obtain some qualitative estimates.
	
	Following this work, many new questions can be asked.  In (Bayesian) optimal experimental design (see, e.g.,~\cite{Alexanderian_BayesianOptExpDesign_2021,ye2013minimizing}), K- and E-optimality seek to maximize the inverse condition number or the minimal eigenvalue of the inverse of the Bayesian covariance matrix or the \gls{fim}, respectively. We examine the same two quantities from a qualitative perspective, and it would be interesting to examine the relation between the two approaches.

	{{As mentioned in Section~\ref{ssec:earlystopping},} our approach suffers from a  drawback that is  {very} typical for all experimental design methods: the sampling of designs requires linearization around a parameter $p_\ast\approx p_0$ close to the ground-truth to build $\tilde \pi$, as demonstrated in \Cref{fig:piLandscapes2}. Several strategies have been developed in classical optimal experimental design literature  to mitigate this drawback of requiring prior knowledge of $p_0$,   summarized in \cite{Alexanderian_BayesianOptExpDesign_2021,huan2024optimalexperimentaldesignformulations} under sequential experimental design, some of which can be directly integrated into our approach. In particular, we see synergies between our approach and the  {myopic} approach consisting of alternating phases of experimental design and parameter inference through gradient-based optimization.} 
	
	Finally, we see potential application of our approach to more recently developed inversion frameworks that rely on a least-squares optimization.  Examples of such frameworks
	can be found in \cite{chen2021solving,dong2023method}, where Gaussian processes or neural networks are incorporated in the inversion process.
	A detailed derivation is  left  for further investigation.

	\section*{Declarations}
	
	\paragraph{Funding}
	K.H. acknowledges support by the German Academic Scholarship Foundation (Studienstiftung des deutschen Volkes), the Marianne-Plehn-Program as well as the Jet Propulsion Laboratory PDRDF 24AW0133. C.K. acknowledges support  from the German Science Foundation, KL566/22-1. Q.L. acknowledges support from DMS-2308440 and DMS-2023239.
	
	\paragraph{Conflict of interest}
	The authors declare no conflicts of interest.
	
	\paragraph{Data Availability}
	Code to generate the examples is available upon request.

	\appendix

		\section{Appendix: Derivation of the formula for $\grd_p u_p(x)$}\label{app:DerivationA(x)}
		
		We derive the formula for the gradient $\grd_p u_p(x)$ of the solution to the Schr\"odinger equation w.r.t. the potential $p$, that we require for the computation of the sampling probabilities.
		
		In the following derivations, all gradients are with respect to $x$, unless specified otherwise. 
		For a fixed measurement location $x\in X$, we can then  define the Lagrange function $\mathcal{L}: \mathcal{A} \times H_0^1(X)\times H_0^1(X)\to \rr$ as
		\[
		\mathcal{L}_{x}(p,u,g) = u(x) +\langle \grd g,\grd u\rangle_{L_2(X)} + \langle g, pu\rangle_{L^2(X)} - \langle g,  {\gamma}\rangle_{H^1(X),H^{-1}(X)} ,
		\]
		where $g$ is the Lagrange multiplier, and $\langle \cdot, \cdot \rangle_{H^1_0(X),H^{-1}(X)}$ denotes the duality bracket in $H^1_0(X)\times H^{-1}(X)$. Using \eqref{forward}, one immediately sees $\mathcal{L}_{x}(p,u_p,g)=u_p(x)$. Therefore, restricted to this solution manifold, the chain rule gives:
		\begin{align*}
			\left.\frac{\partial u_p(x)}{\partial p_j} \right|_{p=\hat p} =& \left.\frac{\partial \mathcal{L}_{x}}{\partial p_j} \right|_{\substack{p=\hat p\\ u = u_{\hat p}}}  + \left.\frac{\partial \mathcal{L}_{x}}{\partial u} \right|_{\substack{p=\hat p\\ u = u_{\hat p}}} \left.\frac{\partial  u_p}{\partial p_j} \right|_{p=\hat p}\,.
		\end{align*}
		This equation holds  for arbitrary $g$, and thus we would like to choose $g = g^{ {(x)}}$ such that $\partial \mathcal{L}_{x}/\partial u = 0$. If so:
		\begin{align*}
			\left.\frac{\partial u_p(x)}{\partial p_j} \right|_{p=\hat p} =& \left.\frac{\partial \mathcal{L}_{x}}{\partial p_j} \right|_{\substack{p=\hat p\\ u = u_{\hat p}}} = \left.\frac{\partial \langle  g^{ {(x)}}, pu\rangle_{L^2(X)}}{\partial p_j}  \right|_{\substack{p=\hat p\\ u = u_{\hat p}}}\\
			=&  \left.\frac{\partial \langle  g^{ {(x)}}, \sum_{k} p_k \phi_k u\rangle_{L^2(X)}}{\partial p_j}  \right|_{\substack{p=\hat p\\ u = u_{\hat p}}}=  \langle  g^{ {(x)}},   \phi_j u_{\hat p}\rangle_{L^2(X)}\,.
		\end{align*}
		It remains to compute $g^{ {(x)}}\in H_0^1(X)$ for which  $\partial \mathcal{L}_{x}(p,u,g^{ {(x)}})/\partial u = 0$. 
		From integration by parts we see
		\begin{align*}
			\partial_u  \mathcal{L}_{x} = &  \partial_u\left[u(x) + \langle\grd g^{ {(x)}}, \grd u\rangle_{L^2(X)}  + \langle g^{ {(x)}}, pu\rangle_{L^2(X)} \right] \\
			=& \partial_u\left[u(x) + \langle -\Delta g^{ {(x)}} ,   u\rangle_{H^{-1}(X),H^1(X)} + \langle p g^{ {(x)}}, u\rangle_{L^2(X)} \right]. 
		\end{align*}
		Setting this to be zero, we have the condition for $g^{ {(x)}}$:
		\begin{align*}
			-\Delta g^{ {(x)}} + p g^{ {(x)}} = -\delta_{x} \text{ on } X, \quad g^{ {(x)}}= 0\text{ on } \partial X\,.
		\end{align*}
		
		\section{Appendix: Proof of \Cref{thm:3matrixproductSketch}}\label{sec:proof3MatrixSketch}
		{For the sake of completeness, we state the proof of \Cref{thm:3matrixproductSketch}, which is a direct adaptation of the proof of Mahoney's result \cite[Theorem 6]{mahoney_2016_RandLA}. 
			\begin{proof}[Proof of \Cref{thm:3matrixproductSketch}]
				Note that $\mathsf X_j$ and  $\hat{\mathsf M}:=  \frac 1 c  \sum_{j=1}^c \mathsf Y_{(k,l)_j} $ are  unbiased estimators of $\sfA ^\top \mathsf W\sfA$ and consider the mean-squared Frobenius error  
				\begin{align*}\mathbb E[\|\sfA ^\top \mathsf W\sfA -\hat{\mathsf M}\|_F^2] = &\sum_{m,n} \operatorname{Var}[\hat{\mathsf M}_{mn}] \leq \frac 1 {c^2} \sum_{m,n}\sum_{j=1}^c \mathbb E[(({\mathsf X_j})_{mn})^2] 
					= \frac 1 c \sum_{k,l} \frac 1{\pi(k,l)} \|\mathsf Y_{k,l}\|^2_F\\
					\leq &\frac 1 {\beta c} \left(\sum_{k,l}  \|\mathsf Y_{k,l}\|_F \right)^2
				\end{align*}
				In order to apply a martingale concentration inequality, one  considers the Frobenius error $\|\sfA^\top \mathsf W \sfA - \hat{\mathsf M}\|_F = :F((k,l)_1,...,(k,l)_c) $ as a function $F$ of the  index pairs $(k,l)_j$ chosen in the realizations of the $\mathsf X_j$. This allows building the Doob martingale  generated by the iterative conditioning on the choice of the $j$-th index pair $(k,l)_j$.  Noting that exchanging  the value of one index pair from $(k,l)_{j_0}$ to  $(k',l')_{j_0}$ while keeping the remaining ones fixed results in controlled  changes in $\hat{\mathsf M} $: 
				$$\|\hat{\mathsf M}- \hat{\mathsf M}'\|_F= \left\|\frac 1 c (\mathsf Y_{(k,l)_{j_0}}- \mathsf Y_{(k',l')_{j_0}})\right\|_F\leq \frac 2 {c\beta} \sum_{k,l} \|\mathsf Y_{k,l}\|_F =:\Delta,$$
				The triangle inequality then gives $$|F((k,l)_1,...,(k,l)_{j_0},...,(k,l)_c)) - F((k,l)_1,...,(k',l')_{j_0},...,(k,l)_c))| \leq \Delta.$$
				Thus, the Hoeffding-Azuma martingale concentration inequality  shows
				$$\mathbb P\left[\|\sfA^\top \mathsf W\sfA - \hat{\mathsf M}\|_F - \frac 1 {\sqrt{\beta c}} \sum_{k,l}\|\mathsf Y_{k,l}\|_F \geq \gamma\right]\leq e^{-\frac {\gamma^2}{2c\Delta^2}},$$
				which yields the high probability bound by a suitable choice of $\gamma = \Delta \sqrt{2c\log(\delta^{-1})}$.
			\end{proof}
		}

	\bibliographystyle{abbrv}
	\bibliography{lit}

	\end{document}